\documentclass[aps,pre,twocolumn,superscriptaddress]{revtex4-1}
\usepackage[dvips]{graphicx}
\usepackage[utf8]{inputenc}
\usepackage{psfrag}
\usepackage{amsthm}
\usepackage{amsmath}
\usepackage{amssymb}

\newtheorem{Theorem}{Theorem}
\newtheorem{Definition}{Definition}

\newtheorem{Corollary}{Corollary}

\newcommand{\sM}{\mathcal{M}}

\newcommand{\abs}[1]{\left\vert #1 \right\vert}
\newcommand{\norm}[1]{\left\Vert #1 \right\Vert}

\newcommand{\R}{{\mathbb R}}  %ams bold
\newcommand{\equivN}{\sim_{\mathcal{N}}}
\newcommand{\equivE}{\sim_{\mathcal{E}}}

\begin{document}

%\preprint{}

\title{Symmetries, Stability, and Control\\ in Nonlinear Systems                                
  and Networks}
  
\author{Giovanni Russo}
\email[]{giovanni.russo2@unina.it}
\affiliation{Department of Systems and Computer Engineering, University of Naples Federico II, Italy. This research was performed while visiting the Nonlinear Systems Laboratory, Massachusetts Institute of Technology}
\author{Jean-Jacques E. Slotine}
\email{jjs@mit.edu}
\affiliation{Nonlinear Systems Laboratory, Massachusetts Institute of Technology, United States}

\date{\today}

\begin{abstract}
  This paper discusses the interplay of symmetries and stability in
  the analysis and control of nonlinear dynamical systems and
  networks. Specifically, it combines standard results on symmetries
  and equivariance with recent convergence analysis tools based on
  nonlinear contraction theory and virtual dynamical systems.  This
  synergy between structural properties (symmetries) and convergence
  properties (contraction) is illustrated in the contexts of network
  motifs arising e.g. in genetic networks, of invariance to
  environmental symmetries, and of imposing different patterns of
  synchrony in a network.
\end{abstract}
%\pacs{*****GR: TO DO}
%\keywords{*****GR: TO DO}
\maketitle

\section{Introduction}

Symmetry is a fundamental topic in many areas of physics and
mathematics~\cite{Gol_Ste_03,Mar_Rad_94,Oli_86}. Many systems in
nature and technology possess some symmetry, which somehow influences
its functionality.  Taking into account such symmetry may simplify the
study of a system of interest. In dynamical systems~\cite{Gol_Ste_03},
symmetry concepts have been used e.g. to explain the onset on
instability in feedback systems~\cite{Meh_05,Meh_07}, to design
observers~\cite{Bon_Mar_Rou_09} and
controllers~\cite{spong,Koo_Mar_97}, and to analyze synchronization
properties and associated symmetry detection
mechanisms~\cite{Pha_Slo_07,Ger_Slo_08}. Typically, the symmetries of a
physical system are preserved in the mathematical tools used to model
it. This is the case for instance of Lagrangian systems, where
can be easily shown that the symmetries of the Lagrangian function
transfer onto the equations of motion, making them invariant under the
same symmetry (see e.g.~\cite{spong} in the context of motion
control).

Our goal in this paper is to develop a theoretical framework to study
the rich interplay between symmetries of the system dynamics and
questions of stability and control. As compared with well-known
fundamental results such as~\cite{Gol_Ste_03}, our approach presents
two key differences. First, in \cite{Gol_Ste_03} the symmetry properties of a
system of interest are studied and its possible final behaviors are
related to these symmetries and to the associated bifurcation
analysis~\cite{Ku:98}.  Our approach is complementary in that it
yields {\it global stability and convergence} results. Second, these
results are further generalized by showing that it is possible to use
\emph{virtual systems} in place of the real systems for performing
convergence analysis. These more general virtual systems may have
symmetries and convergence properties that the real systems do not.

%We first recall standard tools~\cite{Gol_Ste_03} to describe questions
%of symmetry, equivariance, and invariance in dynamical systems. We
%formally show that the final behavior of a dynamical system is
%determined both by its symmetries and by some convergence property of
%its vector field. We then generalize the results linking stability and
%symmetries, and show that some specific behaviors of a system can be
%explained in terms of symmetry and convergence properties of some
%appropriately constructed \emph{virtual} system.

Convergence analysis is based on nonlinear contraction theory
\cite{Loh_Slo_98,Wan_Slo_05}, a viewpoint on incremental stability
which has emerged as a powerful tool in applications ranging from
Lagrangian mechanics to network control. Historically, ideas closely
related to contraction can be traced back to~\cite{Hartmann} and even
to~\cite{Lewis}.  As pointed out in~\cite{Loh_Slo_98}, contraction is
preserved through a large variety of systems combinations, and in
particular it represents a natural tool for the study and design of
synchronization mechanisms~\cite{Wan_Slo_05}.  Here the contraction
theory framework also shows that in fact, for symmetry to play a key
role in convergence analysis and control, it needs not be exhibited by
the physical system itself but only by a much more general {\it
  virtual} system derived from it.  As such, our results provide a
systematic framework extending and generalizing the results
of~\cite{Pha_Slo_07,Ger_Slo_08} in this context.

%In the context of this paper, the contraction theory
%framework also shows that in fact, for symmetry to play a key role in
%convergence analysis and control, it needs not be exhibited by the
%physical system itself but only by a much more general {\it virtual}
%system derived from it.

The paper is organized as follows. After reviewing symmetries and
contraction in Section \ref{sec:math_tools} and Section \ref{sec:symmetries_intro}, some basic results linking
the two notions are described in Section \ref{sec:basic_res_symm_cont}. These results are generalized in
Section \ref{sec:multiple_virtual}, where systems with multiple symmetries are considered. Section \ref{sec:virt_sys_analysis} considerably extends the basic results by showing that the contraction and symmetry conditions on the system of interest can be replaced by weaker conditions on some appropriately constructed virtual system. In Section \ref{sec:from_analys_to_contr}, the approach is applied to the case of systems with
external inputs, with examples detailed in Section \ref{sec:input_scaling}. Using our approach we explain the onset of the so-called fold change detection behavior which is important for biochemical processes. Section \ref{sec:networked_systems} extends our theoretical framework to the study
of interconnected systems or networks, and shows that it can be used
to analyze/control synchronization patterns. Applications are then provided by showing that symmetries and contraction can be controlled so as to generate different synchronization patterns. Quorum sensing networks are also analyzed. Brief concluding remarks are
offered in Section \ref{sec:conclusions}.

\subsubsection*{Notation}
We denote with $\abs{x}$ any vector norm of the vector $x \in \R^n$ and with $\norm{A}$ the induced matrix norm of the real square matrix $A \in \R^{n\times n}$. When needed, we will point out the particular norm being used by means of subscripts: $\abs{\cdot}_i$, $\norm{\cdot}_i$. Given a vector norm on Euclidean space, $\abs{\cdot}$, with its induced matrix norm $\norm{A}$, the associated \emph{matrix measure} $\mu$ is defined as the directional derivative of the matrix norm, that is, $\mu(A) =
\lim_{h \searrow 0} \frac{1}{h} \left(\norm{I+hA}-1\right).$
The matrix measure, also known as \emph{logarithmic norm} was introduced in \cite{dahlquist} and \cite{lozinskii}.
When needed, we will point out the particular matrix measure being used by means of subscripts. Examples of matrix measures are listed in Table
\ref{tab:matrix_measures}. More generally,  matrix measures can be induced  by  weighted vector norms
$\abs{x}_{\Theta,i}= \abs{\Theta x}_i$, with $\Theta$ a constant
invertible matrix and $i=1,2,\infty$. Such measures, denoted with
$\mu_{\Theta,i}$, are linked to the standard measures by $\mu_{\Theta,i}(A)= \mu_i\left(\Theta A \Theta^{-1}\right)$, $\forall i=1,2,\infty$

\begin{table}[th] 
\caption{Standard matrix measures for a real $n\times n$ matrix, $A=[a_{ij}]$. The $i$-th eigenvalue of $A$ is denoted with $\lambda_i(A)$.}
\centering
\label{tab:matrix_measures}
\begin{tabular}{|c| c|} 
\hline
vector norm, $\abs{\cdot}$ & induced matrix measure, $\mu\left(A\right)$\\
\hline
$\abs{x}_1= \sum_{j=1}^n\abs{x_j}$ & $ \mu_1\left(A\right)= \max_{j} \left( a_{jj}+\sum_{i \ne j}  \abs{a_{ij}} \right)$\\
\hline
$\abs{x}_2= \left( \sum_{j=1}^n\abs{x_j}^2\right)^{\frac{1}{2}}$ & $ \mu_2\left(A\right)=\max_{i} \left( \lambda_i\left\{\frac{A+A^\ast}{2}\right\}\right)$\\
\hline
$\abs{x}_\infty= \max_{1 \le j \le n} \abs{x_j}$ & $\mu_{\infty}\left(A\right)= \max_{i} \left( a_{ii}+\sum_{j \ne i} \mid a_{ij}\mid \right)$\\
\hline
\end{tabular}
\end{table}

\section{Mathematical preliminaries}\label{sec:mathematical_preliminaries}

\subsection{Contraction theory tools}\label{sec:math_tools}

The basic result of nonlinear contraction analysis \cite{Loh_Slo_98}
which we shall use in this paper can be stated as follows.

\begin{Theorem}[Contraction]
\label{theorem:contraction}
Consider the $m$-dimensional deterministic system
\begin{equation}
\label{equ:main}
\dot x = f(x,t)
\end{equation}
where $f$ is a smooth nonlinear function.  The system is said to be
{\it contracting} if any two trajectories, starting from different
initial conditions $x_0 = x(t=0)$, converge exponentially to each other. A sufficient
condition for a system to be contracting is the existence of some matrix measure, $\mu$, such that $\exists \lambda > 0$,  $\forall x$,  $\forall t \ge 0$,  $\mu\left(\frac{\partial f(x,t)}{\partial x}\right)  \le  - \lambda$. The scalar $\lambda$ defines the contraction rate of the system.
\end{Theorem}
%
%Note that for linear time-invariant systems, contraction is
%equivalent to strict stability, and, using the Euclidean vector norm,
%$\Theta$ can be chosen as the transformation matrix which
%diagonalizes the system or puts it in Jordan form~\cite{Loh_Slo_98}.

For convenience, in this paper we will also say that a {\it function}
$f(x,t)$ is contracting if the system $\ \dot x= f(x,t) \ $
satisfies the sufficient condition above.  Similarly, we
will then say that the corresponding Jacobian {\it matrix} $\ \frac{\partial
f}{\partial x}(x,t)$ is contracting.

We shall also use the following property of contracting systems,
whose proofs can be found in~\cite{Loh_Slo_98}, \cite{Slo_03}.

\noindent {\bf Hierarchies of contracting systems}\ \ \ Assume that 
the Jacobian of (\ref{equ:main}) is in the form
\begin{equation}\label{eqn:hierarchy_general}
\frac{\partial f}{\partial x}(x,t)\ = \
\left[\begin{array}{*{20}c}
J_{11} & J_{12}\\
0 & J_{22}\\
\end{array}\right]
\end{equation}
corresponding to a hierarchical dynamic structure. The $J_{ii}$ may be
of different dimensions. Then, a sufficient condition for the system
to be contracting is that (i) the Jacobians $J_{11}$, $J_{22}$ are
contracting (possibly with different $\Theta$'s and for different
matrix measures), and (ii) the matrix $J_{12}$ is bounded.

A simple yet powerful extension to nonlinear contraction theory is the concept
of \emph{partial} contraction, which was introduced in~\cite{Wan_Slo_05}. 

\begin{Theorem}[Partial contraction]
\label{theorem:partial_contraction}
Consider a smooth nonlinear $n$-dimensional system of the form $\dot x= f(x,x,t)$ and assume that the auxiliary system $\dot y= f(y,x,t)$
is contracting with respect to $y$.  If a particular solution of
the auxiliary $y$-system verifies a smooth specific property,
then all trajectories of the original $x$-system verify this
property exponentially. The original system is said to be \emph{partially
contracting}.
\end{Theorem}

Indeed, the virtual $y$-system has two
particular solutions, namely $y(t) = x\left(t\right)$ for all $t \ge 0$ and the particular solution with the specific property. Since all trajectories
of the $y$-system converge exponentially to a single trajectory,
this implies that $x\left(t\right)$ verifies the specific property
exponentially.

Using the Euclidean norm, the results in \cite{Wan_Slo_05} are
systematically extended in \cite{Pha_Slo_07} to global exponential
convergence towards some flow-invariant linear subspace, $\sM$, allowing
in particular multiple groups of synchronized elements to co-exist (so
called poly-dynamics, or poly-rhythms).  The dynamics (\ref{equ:main}) is said to be {\emph contracting
towards} $\sM$ if all its trajectories converge towards $\sM$
exponentially. Let $p$ be the dimension of $\sM$ and $V$ be a $(n-p)
\times n$ matrix, whose rows are an orthonormal basis of
$\sM^{\perp}$. The following result is a straightforward
generalization of Theorem 1 in \cite{Pha_Slo_07}:
\begin{Theorem}\label{thm:contraction_invariant_subspace}
If $\mu(V\frac{\partial f}{\partial x}V^T)$ is uniformly negative for
some matrix measure $\mu$ in $\R^{n-p}$, then (\ref{equ:main}) is contracting towards $\sM$. 
\end{Theorem}
Note that if the system is contracting, then trivially it is contracting towards $\sM$
(since entire trajectories of the system are contained in $\sM$).

\subsection{Symmetry of dynamical systems}\label{sec:symmetries_intro}

In this paper, we consider operators acting over the state space of
(\ref{equ:main}). Often such operators are linear, with their effects
on the structure of the solutions specified in terms of a group of
transformations, see e.g. \cite{Gol_Ste_03}. We will use the following
standard definitions.
\begin{Definition}
Let $\Gamma$ be a group of operators acting on $\R^n$. We say that $\gamma \in \Gamma$ is a symmetry of (\ref{equ:main}) if for any solution, $x(t)$, $\gamma x(t)$ is also a solution. Furthermore, if $\gamma x = x$, we say that the solution $x(t)$ is $\gamma$-symmetric.
\end{Definition}
\begin{Definition}
Let $\Gamma$ be a group of operators acting on $\R^n$, and $f:\R^n\times \R^+ \rightarrow \R^n$. The vector field, $f$, is said to be $\gamma$-equivariant if $f(\gamma x,t) = \gamma f(x,t)$, for any $\gamma \in \Gamma$ and $x \in \R^n$.
\end{Definition}
  Thus,  $\gamma$-equivariance in essence means that $\gamma$ ``commutes'' with $f$.
\begin{Definition}
We say that a solution of (\ref{equ:main}) is $h$-symmetric, if there exists some $T>0$ such that $x(t) =\gamma x(t+T)$. The vector field, $f$, is said to be $h$-equivariant if $ f(\gamma x,t)=\gamma f(x,t+T)$.
\end{Definition}

We will refer to $\gamma$ and $h$ as actions.
\subsubsection*{Symmetries, equivariance and invariant subspaces}

We first review the relationship~\cite{Gol_Ste_03} between symmetries,
equivariance, and the existence of flow-invariant linear subspace.

If $f$ is $\gamma$-equivariant, then $\gamma$ is a symmetry of
(\ref{equ:main}). Indeed, letting $y(t) = \gamma x(t)$, we have
$$
\dot y = \gamma \dot x = \gamma f(x,t) = f(\gamma x,t) = f(y,t)
$$
so that $y(t)$ is also a solution of (\ref{equ:main}).

If the operator $\gamma$ is linear, this in turn immediately implies that the subspace
$\sM_{\gamma}=\left\{x \in \R^n:\gamma x = x\right\}$ is
flow-invariant under the dynamics (\ref{equ:main}).  Thus, solutions
having symmetric initial conditions, $x_0 = \gamma x_0$,
preserve that symmetry for any $t\ge 0$. Note that $\sM \neq \emptyset$ since ${\bf 0} \in \sM$. 

In this  paper we assume $\gamma$ to be any linear operator and give some
extensions for nonlinear operators. Therefore, our framework is
somewhat broader than that typically considered in the literature on
symmetries of dynamical systems, where it is generally assumed that
$\gamma$ describes finite groups or compact Lie Groups (see
e.g. \cite{Gol_Ste_03} and references therein).

\subsection{Basic results on symmetries and contraction}\label{sec:basic_res_symm_cont}

Next, we review some results from \cite{Ger_Slo_08} which this paper
shall generalize. These results can be summarized as follows: (i) if the dynamical system of interest is contracting, then $\gamma$ and $h$ symmetries of the vector fields are transferred onto symmetries of trajectories; (ii) if $f$ presents a spatial symmetry $\gamma$, then this property can be
  transferred to the solutions $x(t)$ by only requiring contraction
  towards $\sM_{\gamma}$, rather than contraction of the entire system.

Note that, although the proofs in \cite{Ger_Slo_08} are presented in the context of Euclidean norms, they can be generalized straightforwardly to other norms.
   
\begin{Theorem}\label{thm:symm_contr}
Consider dynamics  (\ref{equ:main}), with $f$ $\gamma$-equivariant.
Assume that $f$ is also contracting, or more generally that $f$ is contracting towards $\sM_{\gamma}$.
Then, any solution of (\ref{equ:main})  converges towards a $\gamma$-symmetric solution.
\end{Theorem}
\proof
The proof is immediate, since any system trajectory tends exponentially towards $\sM_{\gamma}$,
and by definition $\sM_{\gamma}$ is the subspace $x=\gamma x$.\endproof

One interesting interpretation of the above theorem is as
follows. Assume that $f\left(x,t\right)$ in (\ref{equ:main}) is
$\gamma$-equivariant. We know (see Section \ref{sec:symmetries_intro}) that
if $x(t)$ is a solution of (\ref{equ:main}), then so is $\gamma x(t)$,
which implies in particular that the subspace $\sM_{\gamma}$ is
flow-invariant under (\ref{equ:main}). Assume now that $f$ is contracting towards $\sM_{\gamma}$. Then, given arbitrary
initial conditions in $x(t)$, both $x(t)$ and $\gamma x(t)$ will tend
to $\sM_{\gamma}$, and therefore will tend to the {\it same}
trajectory, since by definition $\sM_{\gamma} = \left\{x\in
  \R^n:x=\gamma x\right\}$. Thus, all trajectories initialized within
a group transformation generated by $\gamma$ represent an equivalence
class which will converge to the same trajectory on $\sM_{\gamma}$.

Furthermore, note that adding to the dynamics (\ref{equ:main}) any
term $\dot{x}_\sM(t) \in \sM_\gamma$ preserves contraction to $\sM_\gamma$. By choosing $\dot{x}_\sM(t)$ to
represent a multistable attractor, this property could be used to
spread out or separate solutions corresponding to different
equivalence classes, in a fashion reminiscent of recent work on image
classification~\cite{Mallat}.

A similar transfer of symmetries of the vector field onto symmetries
of $x$ holds for spatio-temporal symmetries. Let $p_{\gamma}$ be the
order of $\gamma$, i.e. $\gamma^{p_\gamma} = \rm{identity}$. The
following result holds:
\begin{Theorem}
  If $f$ is $h$-equivariant and contracting, then $x$ tends to an 
  $h$-symmetric. Furthermore, all the solutions of the
  system tend to a periodic solution of period $p_\gamma T$.
\end{Theorem}
\proof Note first that if $x(t)$ is a solution of (\ref{equ:main}), then
so is $\gamma x(t+T)$, since
$$
\frac{d\gamma x(t+T)}{dt} = \gamma \dot x (t+T) = f (\gamma x(t+T),t)
$$
Since (\ref{equ:main}) is contracting, this implies that
$x(t) \ \rightarrow  \ \gamma x(t+T)\ $ exponentially. By recursion,
$$
x(t) \ \rightarrow \ \gamma^{p_\gamma} x(t+p_\gamma T) = x(t+p_\gamma T)\ \ \ \ \ \ \ {\rm exponentially}
$$
Now exponential convergence of the above implies implies in turn that for any $t \in \left[0,p_\gamma T\right]$, $x(t+np_\gamma
T)$ is a Cauchy sequence. Since $\R^n$ (equipped with either of the weigthed
$1$, $2$ or $\infty$ norms) is a complete space, this shows that the limit
$\lim_{n\rightarrow +\infty} x(t+np_\gamma T)$
does exist, which completes the proof.\endproof

Note that $p_\gamma T$ may actually be an integer multiple of the {\it
  smallest} period of the solutions.

\section{Multiple symmetries and virtual systems}\label{sec:multiple_virtual}

In this Section, we start by extending the results presented above
by considering the case where $f$ presents more than one symmetry. A
further generalization is then given using virtual systems: in this way, our approach is extended to the study of systems which present no symmetries.
 
\subsection{Coexistence of multiple spatial symmetries}\label{sec:multiple_spatial}

In the previous Section, we showed that the symmetries of the vector field of (\ref{equ:main}) are transformed in symmetries of its solutions, $x(t)$, if the system is contracting (towards some linear invariant subspace). We now assume that $f$ is equivariant with respect to a number of $s>1$ actions: the aim of this Section is to provide sufficient conditions determining the final behavior of the system. 

Let: (i)$\sM_i$ be the linear subspace defined by $\gamma_i$ (i.e. $\sM_i = \left\{x : x=\gamma_ix\right\}$); (ii) $\dot x^i = f^i(x^i,t)$ be the dynamics of (\ref{equ:main}) reduced on $\sM_i$; (iii) $\gamma_1,\ldots,\gamma_s$ be the symmetries showed by $f^i$.
\begin{Theorem}\label{thm:multiple_symmetries_1}
Assume that  $\sM_1 \subset \sM_2 \subset \ldots \subset \sM_s$. Then, all the solutions of (\ref{equ:main}) exhibit the symmetry $\gamma_j$ ($1\le j\le s$) if:
(i) (\ref{equ:main}) contracts towards $\sM_s$; (ii) $\forall i=j+1,\ldots, s$, $\dot x^{i} = f^{i}(x^{i},t)$  is contracting towards $\sM_{i-1}$.
\end{Theorem}
\proof
By assumption we know that the sets $\sM_i$ are all linear invariant subspaces. Denote with $\lambda_i$ the contraction rates of $\dot x^{i} = f^{i}(x^{i},t)$ towards $\sM_{i-1}$. Let $a_i(t)$ be solutions of (\ref{equ:main}) such that $a_i\left(0\right) \in \sM_i$, and let $b(t)$ be a solution of (\ref{equ:main}) such that $b\left(0\right) \notin \sM_s$. We have:
$$
\begin{array}{*{20}l}
\abs{b(t)-a_j(t)} =\\
= \abs{b(t) + \sum_{i=j+1}^{s}a_i(t) - \sum_{i=j+1}^{s}a_i(t) - a_j(t)} \le \\
\le \abs{b(t)- a_s(t)} + \sum_{i=j+1}^{s} \abs{a_i(t)- a_{i-1}(t)}
\end{array}
$$
Now, by hypotheses, the dynamics of (\ref{equ:main}) reduced on each of the subspaces $\sM_{i}$ ($i=j+1,\ldots,s$), i.e. $\dot x^{i} = f^i(x^i,t)$, is contracting towards $\sM_{i-1}$. Thus, there exists some $K_i >0$, $i=1,\ldots, j-1$, such that:
$$
\begin{array}{*{20}l}
\abs{b(t)- a_s(t)} \le K_{s+1} e^{-\lambda_{s+1} t}\\
\abs{a_i(t)- a_{i-1}(t)} \le K_{i} e^{-\lambda_{i} t} & i=j+1,\ldots, s
\end{array}
$$
This implies that $\abs{b(t) - a_j(t)} \rightarrow 0$
exponentially. The Theorem is then proved.
\endproof

With the following result, we show that if (\ref{equ:main}) is contracting towards $\sM_\gamma$, then the only symmetries that the vector field, $f$, can eventually exhibit are those defining invariant subspaces strictly contained in $\sM_\gamma$.
\begin{Theorem}\label{thm:multiple_symmetries_2}
Assume that (\ref{equ:main}) exhibits a symmetry, $\gamma$, and that it is contracting towards $\sM_\gamma$. Then there does not exist any other symmetry, $\beta$, such that $\sM_\gamma \cap \sM_\beta = \left\{0\right\}$.
\end{Theorem}
\proof
The proof of this result is straightforward, by contradiction. Assume that there exist two symmetries, $\gamma$ and $\beta$, such that $\sM_\gamma \cap \sM_\beta = \left\{0\right\}$. That is, let 
$d(a,b,t) : =\abs{a(t)-b(t)}$
be the distance between $a(t) \in \sM_\gamma$ and $b(t) \in \sM_\beta$. We have that
$\inf_{a \in \sM_\gamma,\\ b \in \sM_\beta, \\ t \in \R^+} \left\{d(a,b,t)\right\} = D >0$
since both $\sM_\gamma$ and $\sM_\beta$ are invariant subspaces. This in turn implies that $a(t)$ and $b(t)$ cannot globally exponentially converge towards each other. That is, the system is not contracting towards $\sM_\gamma$. This contradicts the hypotheses.
\endproof

With the following result we address the case where the invariant subspaces defined by the symmetries are not strictly contained in each other but intersect:
\begin{Theorem}\label{thm:multiple_symm_intersection}
Assume that $\sM_{\cap}=\cap \sM_i \ne \left\{0\right\}$. Then, all solutions of (\ref{equ:main}) exhibit the symmetry defined by $\sM_{\cap}$ if one of the two conditions holds: (i) $f$ is contracting toward each subspace $\sM_i$; (ii) $f$ is contracting.
\end{Theorem}
\proof 
Let $x_i$, $i=1,\ldots,s$, be solutions of (\ref{equ:main}), such that $x_i(0) \in \sM_i$, and $a(t)$ be a solution of the system such that $a(0) \notin \sM_i$. Now, if $f$ is contracting towards each $\sM_i$, we have, by definition, that there exists $K_i>0$, $\lambda_i >0$, $i=1,\ldots, s$, such that $\abs{a(t)-x_i} \le K_i e^{-\lambda t}$
This, in turn, implies that there exists some $K>0$, $\lambda >0$ such that $\abs{x_i - x_j} \le K e^{-\lambda_i t}, \quad \forall i \ne j$
Now, since $\sM_i$ are flow invariant, we have that $x_i(t) \in \sM_i$, for all $t\ge 0$. Thus, $x_i(t) \rightarrow \sM_{\cap}$, as $t\rightarrow + \infty$, implying that also $a\rightarrow \sM_{\cap}$, as $t \rightarrow + \infty$.

By using similar arguments, it is possible to prove the result under the stronger hypothesis of $f$ being contracting.
\endproof

\subsubsection{Synchronizing networks with chain topologies}\label{sec:chain_top}

As a first application of our results we revisit the problem of finding sufficient conditions for the synchronization of networks having chain topologies. Specifically, we show that Theorem \ref{thm:multiple_symmetries_1} allows to study network synchronization iteratively reducing the dimensionality of the problem.  For the sake of clarity we now consider a simple network of $4$ nodes. While developing the example, we will also introduce introduce an important $\gamma$-symmetry, i.e. permutations. 

\begin{figure}[thbp]
\begin{center}
  \includegraphics[width= 6cm]{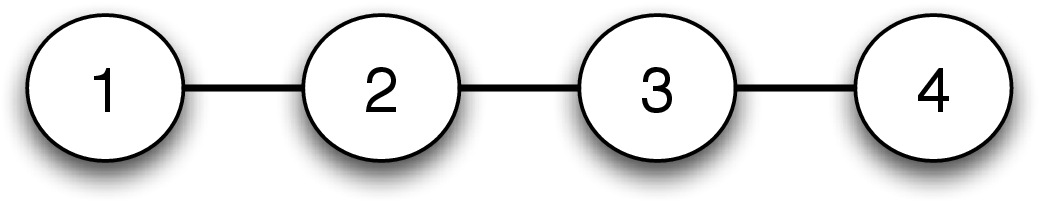}\\
  \includegraphics[width= 6cm]{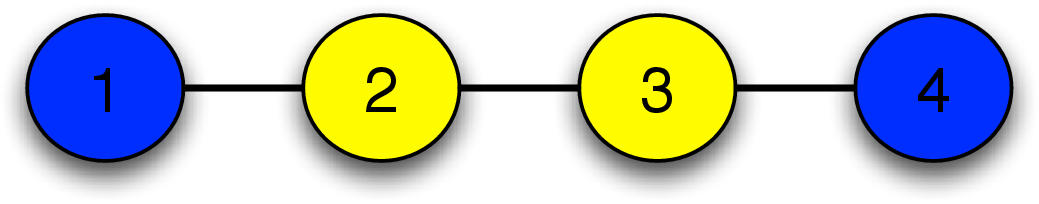}\\
    \includegraphics[width= 3cm]{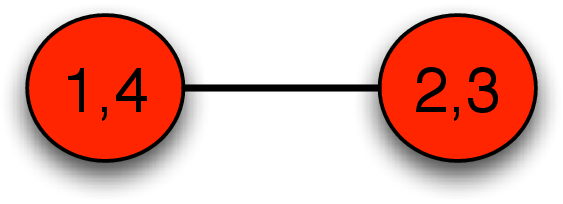}
  \caption{Top panel: the chain topology network of $4$ nodes. Middle panel: poly-synchronous subspace identified by $\sM_2$. Bottom panel: equivalent network and synchronous subspace identified by $\sM_1$.}
  \label{fig:chain_topology}
  \end{center}
\end{figure}

Consider the diffusively coupled network represented in Figure \ref{fig:chain_topology}, whose dynamics are described by:
\begin{equation}\label{eqn:chain_of_4}
\begin{array}{*{20}l}
\dot x_1 = f_1(X) = g(x_1) + h(x_2) - h(x_1)\\
\dot x_2 = f_2(X) = g(x_2) + h(x_1)+h(x_3)-2h(x_2)\\ 
\dot x_3 = f_3(X) = g(x_4) + h(x_2)+h(x_4)-2h(x_3)\\
\dot x_4 = f_4(X) = g(x_4) + h(x_3)-h(x_4)\\
\end{array}
\end{equation}
where: $x_i \in \R^n$, $X=[x_1^T,x_2^T,x_3^T,x_4^T]^T$, all the nodes have the same intrinsic dynamics, $g$ and are coupled by means of the output function, $h$. Now, consider the following action:
\begin{equation}\label{eqn:symmetry_perm_1}
\gamma_2: \quad \left(x_1,x_2,x_3,x_4\right)\rightarrow \left(x_4,x_3,x_2,x_1\right) 
\end{equation}
That is, $\gamma_2$ permutes $x_1$ with $x_4$ and $x_2$ with $x_3$. Let $F(X)=[f_1(X)^T,f_2(X)^T,f_3(X)^T,f_4(X)^T]^T$: it is straightforward to check that $\gamma_2 F(X) = F(\gamma_2 X)$. That is, $F$ is $\gamma_2$-equivariant. This, in turn, implies the existence of the flow invariant subspace
$$
\sM_2 =\left\{X\in \R^{4n}: (x_1,x_2,x_3,x_4) = (x_4,x_3,x_2,x_1)\right\}
$$
Notice that such a subspace corresponds to the poly-synchronous subspace, where nodes $1$ is synchronized to node $4$ and node $2$ is synchronized to node $3$ (synchronous nodes are also pointed out in Figure \ref{fig:chain_topology}). Let $J_2(X)$ be the Jacobian of the network, and 
$$
V_2=\frac{1}{\sqrt{2}}\left[\begin{array}{*{20}c}
-1 & 0 & 0 & 1\\
0 & -1 & 1 & 0\\
\end{array}\right]
$$
be the matrix spanning the null of $\sM_2$ (notice that the rows of $\sM_2$ are orthonormal). All the trajectories of the network globally exponentially converge towards $\sM_2$ if the matrix $V_2J_2(X)V_2^T$ is contracting (see Theorem \ref{thm:symm_contr}). It is straightforward to check that such a matrix is contracting if the function $g(\cdot) - h(\cdot)$ is contracting. Let $x_{1,4},x_{2,3} \in \sM_2$, with $x_{1,4}=x_1=x_4$ and $x_{2,3}=x_2=x_3$; the dynamics of (\ref{eqn:chain_of_4}) reduced on $\sM_2$ is given by:
\begin{equation}\label{eqn:reduced_chain_4}
\begin{array}{*{20}c}
\dot x_{1,4} = g(x_{1,4})+ h(x_{2,3}) - h(x_{1,4})\\
\dot x_{2,3} = g(x_{2,3})+ h(x_{1,4}) - h(x_{2,3})\\
\end{array}
\end{equation}
which corresponds to an \emph{equivalent} $2$-nodes network (see Figure \ref{fig:chain_topology}). It is straightforward to check that the above reduced dynamics is $\gamma_1$-equivariant with respect to the action
$$
\gamma_1 :\quad (x_{1,4},x_{2,3})\rightarrow (x_{2,3},x_{1,4})
$$
Thus, the subspace
$$
M_1 =\left\{ X \in \R^{4n} : (x_{1,4},x_{2,3}) =  (x_{2,3},x_{1,4})\right\}
$$
is a flow invariant subspace. Furthermore, the trajectories of (\ref{eqn:reduced_chain_4}) globally exponentially converge towards $\sM_1$ if $V_1J_1(X)V_1^T$ is contracting, where $V_1 = \frac{1}{\sqrt{2}}[-1 , 1]$ and $J_1(X)$ is the Jacobian of (\ref{eqn:reduced_chain_4}). Now, $V_1J_1V_1^T = \frac{1}{2}(\frac{\partial g}{\partial x_{1,4}}- 2\frac{\partial h}{\partial x_{1,4}}+ \frac{\partial g}{\partial x_{2,3}}- 2\frac{\partial h}{\partial x_{2,3}})$, which is contracting if $g(\cdot) - h(\cdot)$ is contracting.

Thus, using Theorem \ref{thm:multiple_symmetries_1}, we can finally conclude that the network synchronizes if the function
\begin{equation}\label{eqn:synchro_cond}
g(\cdot) - h(\cdot)
\end{equation}
is a contracting function. Furthermore, the synchronization subspace is unique by means of Theorem \ref{thm:multiple_symmetries_2}. 

We remark here that:
\begin{itemize}
\item the dimensionality-reduction methodology presented above can be also extended to the more generic case of chain topologies of length $2^r$, for any integer, $r$;
\item the same methodology can be used to prove synchronization of networks having \emph{hypercube} topologies, as they can be seen as \emph{chains} of \emph{chains}. Hence, the above approach can be used to find condition for the synchronization of \emph{lattices}. Such a topology typically arise from e.g. the discretization of partial differential equations. In this view, our results provide a sufficient condition for the spatially uniform behavior in reaction diffusion $PDE$s, similarly to \cite{Arc_unp};
\item the synchronization condition (\ref{eqn:synchro_cond}) is less stringent than that  obtained by proving contraction of (\ref{eqn:chain_of_4}) towards the synchronization subspace $\sM =\left\{X\in \R^{4n}:x_1=x_2=x_3=x_4\right\}$. 
\end{itemize}

\subsection{Generalizations using virtual systems}\label{sec:virt_sys_analysis}

The results presented in the previous sections link the symmetries of
a dynamical system and contraction. Specifically, they show that if
a system presents a set of $s>1$ symmetries, then the final
behavior is determined by the contraction properties of the vector
field.

In this Section, we extend the previous results and show that in order
for the solutions of (\ref{equ:main}) to exhibit a specific symmetry,
equivariance and contraction of $f$ are not necessarily
needed. Indeed, such a condition can be replaced by a weaker condition: namely, an equivariance
condition on the vector field of some \emph{auxiliary} (or virtual)
system, similar in spirit to Theorem~\ref{theorem:partial_contraction}.

\begin{Theorem}\label{thm:virtual_system_analysis}
Consider the system
\begin{equation}\label{eqn:virtual_system}
\dot y =v(y,x,t)
\end{equation}
where $x(t)$ are the solutions of (\ref{equ:main}) and $v(x,y,t)$ is some smooth function such that:
$$
v(x,x,t)= f(x,t)
$$
The following statements hold:
\begin{itemize}
\item if $v(y,x,t)$ is $\gamma$-equivariant and contracting towards $\sM_{\gamma}$, then any solution of (\ref{equ:main}) converges towards a $\gamma$-symmetric solution;
\item if $v(y,x,t)$ is $h$-equivariant and contracting, then any solution of (\ref{equ:main}) converges towards a $h$-symmetric solution.
\end{itemize}
System (\ref{eqn:virtual_system}) is termed as \emph{virtual system}.
\end{Theorem}
\proof
Indeed, by assumption, all the solutions $y(t)$ of the virtual system globally exponentially converge towards some $h$ ($\gamma$) symmetric solution, say $x(t)$. 
Now, notice that any solution of (\ref{equ:main}), say $a(t)$, is a particular solution of (\ref{eqn:virtual_system}), since $v(x,x,t)=f(x,t)$. This implies that:
$$
\abs{a(t) - x(t)}\rightarrow 0
$$
as $t\rightarrow +\infty$. The result is then proved.
\endproof

Note that
\begin{itemize}
\item Any solution of the virtual system having symmetric initial
  conditions, i.e. $y(0) = \gamma y(0)$, preserves the symmetry for any
  $t > 0$. In particular, if a solution of the {\it real} system has
  initial conditions verifying the symmetry of the {\it virtual} system,
  i.e. $x_0 = \gamma x_0$, then it preserves this symmetry, i.e. $x(t) = \gamma x(t)$ for any $t
  \ge 0$. This is a generalization of the basic result presented in
  Section \ref{sec:symmetries_intro}.
\item  Theorem \ref{thm:virtual_system_analysis} can be straightforwardly extended to the case where the virtual system presents a set of $s>1$ spatial symmetries. Analogous results to theorems \ref{thm:multiple_symmetries_1} - \ref{thm:multiple_symm_intersection} can be easily proven.
\end{itemize}

%To briefly illustrate the above result, consider the following simple dynamical system (see Section \ref{sec:chem_synchro} for another application of the result):
%\begin{equation}\label{eqn:virt_example_1}
%\dot x = - a x + b(x) x^4
%\end{equation}
%where $a>0$ and the function $b(x)$ is such that $b(x)x\le0$ for any $x$. We assume that the above system presents no symmetries. Now, consider the following virtual system
%\begin{equation}\label{eqn:virt_example_2}
%\dot  y = -ay+b(x)xy^3
%\end{equation}
%Notice that the virtual system is $\gamma$-equivariant (for any choice of the function $b(x)$) with respect to the action
%$$
%\gamma : \quad y \rightarrow -y
%$$
%Furthermore, the Jacobian of (\ref{eqn:virt_example_2}) is given by:
%$$
%J= -a + 3b(x)xy^2
%$$
%Therefore, the virtual system (\ref{eqn:virt_example_2}) is contracting for any $x$. This, in turn implies that all of its solutions globally exponentially converge to $\sM_\gamma = \left\{ y =0 \right\}$. Thus, $x(t)\rightarrow 0$, as $t\rightarrow + \infty$.

\subsubsection*{A discussion on symmetries of virtual systems}

Let us briefly discuss some of the main features of our
results involving the use of virtual systems.

We showed that a given dynamical system of interest can exhibit some
symmetric final behavior even if the corresponding vector field
is not equivariant and/or contracting. Indeed, a sufficient condition
for a system to exhibit a symmetric final behavior is the symmetry of
the vector field of some appropriately constructed virtual system. Of
course, an interesting general question is that of identifying a
virtual system explaining the final behavior of a real system, an
aspect is reminiscent of the process of identifying a Lyapunov
function in stability analysis.

The idea of relating behaviors of real systems
using a symmetric \emph{virtual} system, possibly of different
dimension, presents analogies with the concept of \emph{supersymmetry}
in particle physics (see e.g. \cite{Fre_88,Ryd_96} and
references therein). The motivation beyond the concept of
supersymmetry is that non-symmetric transformations of an object (the
real system in our framework) in a finite dimensional space, may be
explained by a symmetric transformation of another, possibly
higher-dimensional object (the virtual system in our framework).

Finally, we remark here that all the results presented above 
can be straightforwardly extended to address the problem of designing control strategies guaranteeing convergence of a system of interest onto some desired trajectory. Intuitively, the idea is that the control input has to: i) generate some desired symmetry for the vector field (and hence some desired invariant subspace defining the system's final behavior); ii) \emph{drive} all the trajectories towards the invariant subspace imposing contraction.

\section{System with inputs}\label{sec:from_analys_to_contr}

In the above sections, we presented some results that can be used to
analyze the final behavior of a system of interest. The main idea
beyond such results is the use of contraction to study convergence of
trajectories towards some invariant subspace. In turn, such a subspace
is defined by some \emph{structural} property of the vector field,
namely a symmetry. 

We now generalize this result further, and show that it is possible to
determine a direct relation between the trajectories of a system when
forced by different classes of inputs. The results presented in this Section are also based on the concept of virtual system. Indeed, while the forced systems of interest considered here are not equivariant and/or contracting, we will show that it is
possible to construct a symmetric and contracting virtual system which
allows us to relate the final behavior of the two systems.

Consider a system described by:
\begin{equation}\label{eqn:system_symmetry}
\dot x = f(x,u(t),t)
\end{equation}
The following result holds:
\begin{Theorem}\label{thm:virtual_system_diff_symmetries}
Assume that (\ref{eqn:system_symmetry}) is contracting with respect to $x$, uniformly in $u(t)$, and that there exist some linear transformations $\gamma_i$, $\rho_i$, $i \ge 1$, such that:
$
\gamma_i f(x, u(t),t) = f(\gamma_i x, \rho_i  u(t),t)
$
 Let $x_i(t)$ be solutions of (\ref{eqn:system_symmetry}) when forced by $u(t) = u_i(t)$, i.e. $\dot x_i = f(x_i, u_i(t),t), \quad x_i(t=0) = x_{0,i}$.
Then, for any $u_i(t)$, $u_j(t)$ such that $\rho_i u_i(t) = \rho_j u_j(t)$
$$
\abs{\gamma_ix_i -\gamma_jx_j}\rightarrow 0
$$
as $t\rightarrow + \infty$. Moreover, let $x_i^k$ and $x_j^k$ the $k$-th component of $x_i$ and $x_j$ respectively and $\gamma_i^k$, ($\gamma_j^k$) be the $k$-th component of $\gamma_i$ ($\gamma_j$). If 
$
\gamma_i^k x_{0,i}^k = \gamma_j^k x_{0,j}^k
$.
then $\gamma_i^k x_i^k(t) = \gamma_j^k x_j^k(t)$, for any $t \ge 0$.
\end{Theorem}

The second statement of the above Theorem implies that if the system when forced by two different inputs starts with certain symmetries, then the symmetries are preserved.
\proof
Indeed, let $u_v =\rho_iu_i = \rho_ju_j$ and consider the following virtual system:
\begin{equation}\label{eqn:virt_sys_proof_1}
\dot y = f(y,u_v,t)
\end{equation}
Notice that, for any $i$, $j$, $\gamma_i x_i$ and $\gamma_jx_j$ are particular solutions of such a system. Indeed:
$$
\begin{array}{*{20}l}
\gamma_i \dot x_i = \gamma_if(x_i, u_i,t) = f(\gamma_ix_i,\rho_i u_i,t) = f(\gamma_ix_i,u_v,t)\\
\gamma_j \dot x_j = \gamma_jf(x_j, u_j,t) = f(\gamma_jx_j,\rho_j u_j,t) = f(\gamma_jx_j,u_v,t)\\
\end{array}
$$
Now, since $f(x,u_v,t)$ is contracting by hypotheses, we have, for any $i$, $j$, there exists some $C$ such that:
$$
\abs{\gamma_ix_i - \gamma_jx_j}\le C \abs{\gamma_ix_{0,i}-\gamma_j x_{0,j}}e^{-\lambda t}, \quad \lambda >0
$$
This proves the first part of the result. To conclude the proof it suffices to notice that exponential convergence of $\abs{\gamma_ix_i - \gamma_jx_j}$ to $0$ implies that all of its components exponentially converge to $0$. In particular, this implies that there exists some $C_k$, $\lambda_k$ such that:
$$
\abs{\gamma_i^kx_i^k - \gamma_j^kx_j^k} \le C_k \abs{\gamma_i^kx_{0,i}^k-\gamma_j^k x_{0,j}^k}e^{-\lambda t_k}, \quad \lambda_k >0
$$
Since $\abs{\gamma_i^kx_{0,i}^k-\gamma_j^k x_{0,j}^k} = 0$ by hypotheses, we have that  $\abs{\gamma_i^kx_{i}^k(t)-\gamma_j^k x_{j}^k(t)}$, $\forall t \ge 0$
\endproof

Theorem \ref{thm:virtual_system_diff_symmetries} can be extended by
replacing the linear operators $\gamma_i$, $\rho_i$ by more general
nonlinear transformations acting on the system
\begin{equation}\label{eqn:system_symmetry_feedback}
\dot x = f(x,u(x,t),t)
\end{equation}
The transformations considered  are smooth
nonlinear functions of the state and of time,
$
\gamma=\gamma(x,t)$, $\rho=\rho(u(x,t),x,t)$
Following the same
arguments as in
Theorem~\ref{thm:virtual_system_diff_symmetries}, it is
then straightforward to show,
\begin{Theorem}\label{thm:virtual_system_diff_symmetries_nl}
Assume that (\ref{eqn:system_symmetry_feedback}) is contracting uniformly in $u(x,t)$ and that there exist some $\gamma_i(x,t)$, $\rho_i(u(x,t),x,t)$, $i \ge 1$, such that:
$$
\frac{\partial\gamma_i}{\partial x} f(x, u(x,t),t) = f(\gamma_i(x), \rho_i(u(x,t),x,t),t)
$$
 Let $x_i(t)$ be solutions of (\ref{eqn:system_symmetry_feedback}) when forced by $u(t) = u_i(x_i,t)$, i.e. $\dot x_i = f(x_i, u_i(x_i,t),t)$, $x_i(t=0) = x_{0,i}$
Then, for any $u_i(x_i,t)$, $u_j(x_j,t)$ such that $\rho_i(u(x_i,t),x_i,t) = \rho_j(u_j(x_j,t), x_j,t)$
$$
\abs{\gamma_i(x_i) -\gamma_j(x_j)}\rightarrow 0
$$
as $t\rightarrow + \infty$. Moreover, let $x_i^k$ and $x_j^k$ the $k$-th component of $x_i$ and $x_j$ respectively and $\gamma_i^k$, ($\gamma_j^k$) be the $k$-th component of $\gamma_i$ ($\gamma_j$). If
$$
\gamma_i^k \left(x_{0,i}^k\right) = \gamma_j^k \left(x_{0,j}^k\right)
$$
then $\gamma_i^k\left( x_i(t)^k\right) = \gamma_j^k \left(x_j(t)^k\right)$, for any $t \ge 0$.
\end{Theorem}

\proof
The proof follows exactly the same steps as those used to prove
Theorem \ref{thm:virtual_system_diff_symmetries}, with $u_v$ in
virtual system (\ref{eqn:virt_sys_proof_1}) now being chosen as
$u_v=\rho_i(u(x_i,t),x_i,t) = \rho_j(u_j(x_j,t), x_j,t)$.
\endproof

%*****GR comment: I included the following remarks to answer to the points we discussed on skype

We close this Section by pointing out some features of the above two theorems. 
\begin{itemize}
\item the proofs of both Theorem \ref{thm:virtual_system_diff_symmetries} and Theorem \ref{thm:virtual_system_diff_symmetries_nl} are based on the proof of contraction of some appropriately constructed virtual system of the form (\ref{eqn:virt_sys_proof_1}). We now show that, if some hypotheses are made on $\gamma_i$'s, then the contraction condition can be weakened. Specifically, assume that all the intersection of the subspaces defined by $\gamma_i$, $\sM_{i}$, is nonempty. Then, it is straightforward to check that $\abs{\gamma_i x_i -\gamma_j x_j}\rightarrow 0$ if: (i) $f$ is contracting towards each $\sM_{i}$, or (ii) contracting towards $\sM_{\cap}$. Notice that, since our results make use of symmetries of virtual systems, they extend those in \cite{Ger_Slo_08};
\item Analogously, Theorem \ref{thm:virtual_system_diff_symmetries} and Theorem \ref{thm:virtual_system_diff_symmetries_nl} can also be extended to study the case where the input $u_i$ \emph{selects} one specific symmetry $\gamma_i$. Indeed, let $u_v =\rho_i u_i$. In this case, it can be shown that symmetry $\gamma_i$ is shown by the solutions of (\ref{eqn:system_symmetry}) if $f(x,u_v,t)$ is contracting towards $\sM_{i}$.
\item A particularly interesting case for Theorem \ref{thm:virtual_system_diff_symmetries} is when some of the components of $x_{0,i}$ and $x_{0,j}$ are the same and the actions $\gamma_i$ and $\gamma_j$ leave such components unchanged. That is, in view of the notations above $\gamma_i^k=id=\gamma_j^k$ and $x_{0,i}^k = x_{0,j}^k$. Indeed, in this case Theorem \ref{thm:virtual_system_diff_symmetries} implies that 
$
x_i^k(t) = x_j^k(t), 
$ $\forall t \ge 0$
That is, the $k$-th components of the trajectories of (\ref{eqn:system_symmetry}) have identical temporal evolutions even if forced by different inputs. A similar result holds for Theorem \ref{thm:virtual_system_diff_symmetries_nl} .This consequence of the above two results is used in Section \ref{sec:input_scaling}.
\end{itemize}

\section{An example: invariance under input scaling}\label{sec:input_scaling}

In a series of recent papers, input-output properties of some cellular signaling biochemical systems have been analyzed \cite{Goe_09,Alon_nd,Goe_Kir_09,Coh_Coh_Sig_Lir_Alo_09}. Such studies point out that many sensory systems show the property of having their output invariant under input scaling, which can be formally defined as follows:
\begin{Definition}\label{defn:FCD}
Let $x_i(t)$, $x_j(t)$ be solutions of (\ref{eqn:system_symmetry}) with initial conditions $x_0=x_{i}(0)=x_j(0)$, when $u(t) = \chi_i(t)$ and $u(t)=\chi_j(t)$, respectively. System (\ref{eqn:system_symmetry}) is invariant under input scaling if $x_i(t) = x_j(t)$ for any $\chi_i(t)$, $\chi_j(t)$ such that $\chi_j(t) = F(t) \chi_i(t)$, with $F(t) > 0$.
\end{Definition}

Invariance under input scaling with {\it constant} $F(t) = F$ has been
recently studied in transcription networks by Alon and his
co-authors~\cite{Goe_09,Coh_Coh_Sig_Lir_Alo_09,Alon_nd}. In such
papers, the authors focus on the study of transcriptional networks
subject to step-inputs. In this case, the invariance under input
scaling is called \emph{fold-change detection behavior} (FCD), as the
output of the system depends only on fold changes in input and not on
its absolute level. For example, if the input to the system is a step
function from $1$ to $2$, then its output is the same as if the step
was increased from $2$ to $4$.

%In \cite{Alon_nd} it is shown that FCD is necessary and sufficient to make sensory searches in which an organism moves through a spatial sensory field invariant to the amplitude of the field. This feature is of fundamental importance in e.g. vision, \cite{Alon_nd}. In fact, as shown in \cite{Kee_Sne_09}, the reflectance of an object, say $R(r)$, is multiplied by the ambient light, $I$, to provide the contrast filed sensed by the eye. The eyes make spatial searches by means of rapid movements (fixational eye movements) several times per second, which scan the visual field. FCD in such a system, would allow visual searches to be independent on the strength of the ambient light. Recent studies suggest that spatial visual searches, where the eyes search for a specific object, are insensitive to ambient-light levels (across several order of magnitude), see e.g. \cite{Wal_Har_Bar_06}.

This section uses this paper's results to analyze the associated mathematical
models, arising from protein signal-transduction systems and bacterial
chemotaxis, and in particular it revisits the recent work~\cite{Alon_nd}
from this point of view.  It also shows how these results could, for
instance, suggest a mechanism for stable quorum sensing in bacterial
chemotaxis, thus combining symmetries in cell interactions (quorum
sensing) with invariance to input scaling (fold change detection).

\subsection{Gene regulation}\label{sec:gene_reg}

\begin{figure}[thbp]
\begin{center}
  \includegraphics[width= 5cm]{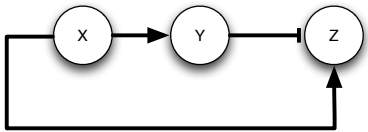}
  \caption{A schematic representation of the I1-FFL}
  \label{fig:I1FFL}
  \end{center}
\end{figure}

This first example considers a pattern (network motif) arising in gene regulation networks, the \emph{Type $1$ Incoherent Feed-Forward Loop} (I1-FFL)  \cite{Eic_04}, \cite{Mil_She_Itz_Kas_Chk_Alo_02}.
The \emph{I1-FFL} is one of the most common network  motifs in gene regulation networks (see also Section \ref{sec:applications}). As shown in Figure \ref{fig:I1FFL}, it consists of an activator, $X$, which controls a target gene, $Z$, and activate a repressor of the same gene, $Y$ (which can be thought of as the output of the system). It has been recently shown that such a network motif can generate a temporal pulse of $Z$ response, accelerate the response time of $Z$ and act as a band-pass amplitude filter, see e.g. \cite{Man_06}, \cite{Kim_08}.

In \cite{Goe_09} it has also been shown by using a dimensionless analysis that for a certain range of biochemical parameters, the I1-FFL can exhibit invariance under step-input scaling (i.e. FCD).

\subsubsection{A basic model}\label{sec:basic_model}

In \cite{Goe_09}, it was shown that a minimal circuit which achieves FCD is the I1-FFL, with the activator in linear regime and the repressor saturating the promoter of the target gene, $Z$.  The model in \cite{Goe_09} is of the form
\begin{equation}\label{eqn:fold_basic_model}
\begin{array}{*{20}l}
\dot Y = -\alpha_1Y+ \chi(t)\\
\dot Z = \beta_2\frac{\chi(t)}{Y}-\alpha_2Z
\end{array}
\end{equation}
where $\alpha_1$, $\alpha_2$, $\beta_2$ are biochemical (positive) parameters and $\chi(t)$ is the input to the system (which can be approximated by the concentration of $X$). It was also shown that the dimensionless model
$$
\begin{array}{*{20}c}
\frac{dy}{d\tau} = F- y\\
r \frac{dz}{d\tau} = \frac{F}{y}-z
\end{array}
$$
with:
$$
\begin{array}{*{20}l}
y = \frac{Y\alpha_1}{\beta_1\chi_{\min}} & Z = \frac{Z}{\beta_2\alpha_1/\beta_1\alpha_2}\\
F= \frac{\chi(t)}{\chi_{\min}} & \tau = \alpha_1 t
\end{array}
$$
exhibits invariance under input scaling.  Later we will also consider more detailed mathematical models that in \cite{Goe_09} have been analyzed numerically. In Section \ref{sec:applications}, we will also analyze other important network motifs under a slightly different viewpoint, i.e. by considering each of the \emph{species} composing the motif as \emph{nodes} of an interconnected systems.

In this Section, we show invariance under input scaling for system (\ref{eqn:fold_basic_model}) for any input, $\chi_i(t)$, $\chi_j(t)$, such that
$
\frac{\chi_i(t)}{\chi_{\min,i}} = \frac{\chi_j(t)}{\chi_{\min,j}} = F(t)
$
In the above expressions $\chi_{\min,i}$ and $\chi_{\min,j}$ denote the basal level of the inputs $\chi_i(t)$ and $\chi_j(t)$ respectively. Such levels are assumed to be nonzero. Notice that the above class of inputs is wider that the one used in Definition \ref{defn:FCD}.

Theorem \ref{thm:virtual_system_diff_symmetries} is now used to prove invariance under input scaling for (\ref{eqn:fold_basic_model}). That is, we show that invariance under input scaling is a consequence of the existence of a symmetric and contracting virtual system in the spirit of Theorem \ref{thm:virtual_system_diff_symmetries}.

In what follows, we will denote with $x_i=(Y_i,Z_i)^T$ and $x_j=(Y_j,Z_j)^T$ the solutions of (\ref{eqn:fold_basic_model}), when $\chi(t) =\chi_i(t)$ and $\chi(t) =\chi_j(t)$, respectively. We assume that $Z_i(0) = Z_j(0)$. In terms of the notation introduced in Theorem \ref{thm:virtual_system_diff_symmetries}, we have $u(t) = \chi(t)$ and:
$$
f(x,u(t)) = \left(\begin{array}{*{20}c}
  -\alpha_1Y+\chi(t)\\
  \beta_2\frac{\chi(t)}{Y}-\alpha_2Z\\
\end{array}\right)
$$
Now, define the following actions:
\begin{equation}\label{eqn:actions}
\begin{array}{*{20}c}
\gamma_i = 
\left(\begin{array}{*{20}c}
Y\\
Z\\
\end{array}\right) \rightarrow
\left(\begin{array}{*{20}c}
\frac{Y}{\chi_{\min,i}}\\
Z\\
\end{array}\right), & \rho_i : \chi(t) \rightarrow \frac{\chi(t)}{\chi_{\min,i}}=F(t)
\end{array}
\end{equation}
It is straightforward to check that:
\begin{itemize}
\item $f(x,u(t))$ is contracting uniformly in $u(t)$;
\item  $\gamma_if(x_i,\chi_i(t)) = f(\gamma_ix_i,\rho_i\chi_i(t))$
\end{itemize}
Now, Theorem \ref{thm:virtual_system_diff_symmetries} implies that for any input $\chi_i(t)$, $\chi_j(t)$ such that $\rho_i\chi_i(t) = \rho_j\chi_j(t)$, $\gamma_ix_i$ and $\gamma_jx_j$ globally exponentially converge towards each other. That is:
\begin{equation}\label{eqn:fold_convergence}
\abs{\gamma_ix_i-\gamma_jx_j} =\abs{\left(\begin{array}{*{20}c}
\frac{Y_i}{\chi_{\min,i}} -\frac{Y_j}{\chi_{\min,j}}\\
Z_i - Z_j
\end{array}\right)} \rightarrow 0
\end{equation}
for any $\chi_i$, $\chi_j$ such that:
\begin{equation}\label{eqn:class_input}
\frac{\chi_i(t)}{\chi_{\min,i}} = \frac{\chi_j(t)}{\chi_{\min,j}} =F(t)
\end{equation}
Now, (\ref{eqn:fold_convergence}) implies that $\abs{Z_i -Z_j} \rightarrow 0$ exponentially. Since the initial conditions of $Z_i$ and $Z_j$ are the same we have that $Z_i(t) = Z_j(t)$ for any $t \ge 0$. That is, the system exhibits invariance under input scaling.

\subsection{A model from chemotaxis}\label{sec:chemo}

In bacterial chemotaxis, bacteria walk through a chemo-attractant
field, say $u(t,r)$ (r denotes two dimensional space vector). Along
their walk, bacteria sense the concentration of $u$ at their position
and compute the tumbling rate (rate of changes of the direction) so as
to move towards the direction where the gradient increases, see
e.g. \cite{Cel_Ver_10}. Typically, the input field is provided by
means of a source of attractant which diffuses in the medium with
bacteria accumulating in the neighborhood of the source. In this case,
the information on the position of the source is encoded only in the
shape of the field and not in its strength. Therefore, it is
reasonable for bacteria to evolve a search pattern which is dependent
only on the shape of the field and not on its strength, i.e., a search
pattern which is invariant under input scaling~\cite{Alon_nd}.
Specifically, consider the following model~\cite{Alon_nd} adapted from
the chemotaxis model of \cite{Tu_Shi_Ber_08}:
\begin{equation}\label{eqn:alon_second_example_1}
\begin{array}{*{20}l}
\ \dot x  = x f(y)\\
\epsilon \dot y  = \phi\left(\frac{u}{x}\right) - y\\
\end{array}
\end{equation}
where $u>0$ is an increasing step-input to the system, representing the ligand concentration, and $y>0$, the output of the system, represents the average kinase activity. The quantity
$x>0$ is an internal variable. We assume the function $\phi$ to be: i) a decreasing function in $x$ with bounded partial derivative $\partial \phi/\partial x$; ii) an increasing function of $u/x$, with derivative  $\phi'=\partial \phi/\partial (u/x) \le b$, $b>0$. Note  that the above model becomes the one used in~\cite{Alon_nd}, when $\phi(u/x)=u/x$. Such a model is obtained assuming $x$ is sufficiently large, with the term $u/x$ actually a simplification of
a  term of the form $u/(x+\eta)$, with $0< \eta \ll x$. The positive constant $\epsilon$ is typically
small, so as to represent a separation of time-scales.

Assume  as in \cite{Tu_Shi_Ber_08} that $f(1) =0$ and that $f(y)$ is strictly increasing with $y$.
Obviously, (\ref{eqn:alon_second_example_1}) verifies the symmetry conditions of Theorem~\ref{thm:virtual_system_diff_symmetries} with:
\begin{equation}\label{eqn:chem_symmetry}
\begin{array}{*{20}c}
\gamma_i :(x_i,y_i)\rightarrow \left(\frac{x_i}{\bar u_i}, y_i\right) & \ \ \ \ \ \rho_i : u_i \rightarrow \frac{u_i}{\bar u_i}
\end{array}
\end{equation}
where $\bar u_i$ denotes the initial (lower) value, at time $t=0$, of the step function. As in the previous Section we assume that $y_i(0) = y_j(0)$. Now, by means of Theorem \ref{thm:virtual_system_diff_symmetries}, we can conclude that, if the system is contracting, 
 $ y_i(t) = y_j(t)$, $\forall t \ge 0$, for any input such that $\frac{u_i}{\bar u_i} = \frac{u_j}{\bar u_j} = F$.

Let us derive a condition for (\ref{eqn:alon_second_example_1}) to be
contracting, which will give conditions on the dynamics and inputs of
(\ref{eqn:alon_second_example_1}) ensuring invariance under input scaling. Model (\ref{eqn:alon_second_example_1}) can be recast as
\begin{equation}\label{eqn:dotdot_y}
\ddot y + \frac{1}{\epsilon}\dot y - \frac{1}{\epsilon}\frac{\partial \phi}{\partial x} x f(y) =0
\end{equation}
As in \cite{Tu_Shi_Ber_08,Alon_nd}, choose $f(y)=y-1$ for simplicity, so that
(\ref{eqn:dotdot_y}) becomes
\begin{equation}\label{eqn:dotdot_y_simplified}
\ddot y + \frac{1}{\epsilon}\dot y -  \frac{1}{\epsilon}\frac{\partial \phi}{\partial x} x (y -1) =0
\end{equation}
The above dynamics is similar to a mechanical mass-spring-damper system with
a time-varying spring,
$\ddot r + 2\eta \omega \dot r + \omega ^2 r =0$
with $2\eta\omega = \frac{1}{\epsilon}$ and $\omega^2 = -
\frac{1}{\epsilon}\frac{\partial \phi}{\partial x}x$.  Now, as shown in \cite{Loh_Slo_unp} such a
dynamics is contracting if $\eta>\frac{1}{\sqrt 2}$. Thus, it immediately follows that (\ref{eqn:dotdot_y_simplified}) is contracting if:
\begin{equation}\label{eqn:cond_contr_chem_gen}
\frac{\partial \phi}{\partial x}x > - \frac{1}{2\epsilon}
\end{equation}
Hence, contraction is attained if $\phi ' \left(-\frac{u}{x^2}\right) x > -\frac{1}{2\epsilon}$
That is, a sufficient condition for (\ref{eqn:dotdot_y_simplified}) to be contracting is
\begin{equation}\label{eqn:cond_contr_chem_gen_2}
x > 2 \epsilon u b
\end{equation}

Notice that, in the case where $\phi\left(u/x\right) =u/x$,  (\ref{eqn:cond_contr_chem_gen_2}) simply becomes:
\begin{equation}\label{eqn:contr_cond_chemotaxis}
x > 2 \epsilon u
\end{equation}
The above inequality implies that, in this case, the system is
contracting (and hence exhibits invariance under input scaling) if the level
of $x$ is sufficiently high (which is true by hypotheses)
and its dynamics is sufficiently slow ($\epsilon$ small) with respect to the dynamics of $y$. Also, given $\varepsilon<\frac{1}{2}$ and a constant $u$, if
contraction condition (\ref{eqn:contr_cond_chemotaxis}) is verified at
$t=0$ with initial conditions embedded in a ball contained in the contraction region (\ref{eqn:contr_cond_chemotaxis}), it remains verified for any $t \ge 0$.

Finally, note that the results of this section, and indeed of the
original~\cite{Goe_09,Alon_nd}, are closely related to the idea, first
introduced in~\cite{Pha_Slo_07} and further studied in
\cite{Ger_Slo_08}, of detecting a symmetry (here, in the environment)
by using a dynamic system having the same symmetry.

\section{Analysis and control of interconnected systems}\label{sec:networked_systems}

The results presented in the previous Section indicate that there exists a direct link between symmetries of a (virtual) vector field and of its solutions, if the system is contracting (or it is made contracting by some control input). 

The aim of this Section is that of using the above results to analyze and control the (poly-) synchronous behavior of $N>1$ (possibly heterogeneous) interconnected systems (also termed as networks in what follows). Such a behavior has been recently reported in ecological systems, networks characterized by strong community structure and in bipartite networks consisting of two groups (see e.g. \cite{Mon_Kur_Bla_04,Oh_Rho_Hon_Kah_05,Sor_Ott_07}). For interconnected systems, symmetries are essentially defined by the nodes' dynamics, the topology of the network and by the particular choice of the coupling functions.

As a first step we will introduce the notion of interconnected system that will be used in the following, see also \cite{Gol_Ste_06}. 

\subsection{Definitions}
In our framework the phase space of the $i$-th node (or cell, or neuron) is denoted with $P_i$, while its state at time $t$ is denoted with $x_i(t)$. Notice that $P_i$ could in general be a manifold. Each node has an intrinsic dynamics, which is affected by the state of some other nodes (i.e. the neighbors of $i$) by means of some coupling function. Those interactions will be represented by means of directed graphs. In such a graph the nodes having the same internal dynamics will be represented with the same symbol. Analogously, heterogeneous coupling functions can be taken into account: identical functions will be denoted by the same symbol.

This is formalized with the following:
\begin{Definition}
An interconnected system consists of: {\it (i)} a set of nodes $\mathcal{N}= \left\{1,\ldots, N\right\}$; {\it (ii)} an equivalence relation, $\equivN$ on $\mathcal{N}$; {\it (iii)} a finite set, $\mathcal{E}$, of edges (arrows); {\it (iv)} an equivalence relation, $\equivE$ on $\mathcal{E}$; {\it (v)} the maps $\mathcal{H}: \mathcal{E} \rightarrow \mathcal{N}$ and $\mathcal{T}: \mathcal{E} \rightarrow \mathcal{N}$ such that: for $e \in \mathcal{E}$, we have $\mathcal{H}(e)$ is the head of the arrow and  $\mathcal{T}(e)$ the tail of the arrow; {\it (vi)} equivalent arrows have equivalent tails and heads. That is, if $e_1,e_2 \in \mathcal{E}$ and $e_1 \equivE e_2$, then $\mathcal{H}(e_1) \equivN \mathcal{H}(e_2)$ and $\mathcal{T}(e_1) \equivN \mathcal{T}(e_2)$.
\end{Definition}

We say that an edge $e \in \mathcal{E}$ is an \emph{input edge} to a node, say $i$, if 
$\mathcal{H}\left(e\right)=i$. The set of  input edges to node $i$ is termed as \emph{input set} and denoted by $\mathcal{I}\left(i\right)$. We also say that two nodes, say $c$ and $d$, are input-equivalent if there exists an arrow type preserving bijection, $\beta : \mathcal{I}(c)\rightarrow \mathcal{I}(d)$.

%The following set of edges, defining an important equivalence relation, is associated to each node, $i$.
%\begin{Definition}
%For any $i \in \mathcal{N}$, the input set of $i$ is defined as
%$$
%\mathcal{I}(i) = \left\{e \in \mathcal{E} : \mathcal{H}(e) = i \right\}
%$$
%Any element of $\mathcal{I}(i)$ is termed as \emph{input edge (or arrow)} of $i$.
%\end{Definition}
%\begin{Definition}\label{def:input_equivalence}
%The relation $\equivI$ (input equivalence) on $\sN$ is defined by $c \equivI d$ if and only if there exists an arrow type preserving bijection $\beta: \sI(c) \rightarrow \sI(d)$.
%\end{Definition}
Finally, our set-up is completed by defining the dynamics of an interconnected system as follows:
\begin{Definition}\label{defn:interconnected_system}
The dynamical system
\begin{equation}\label{eqn:coupled_sys}
\dot X = F(X,t)
\end{equation}
defines an interconnected system if its phase space is defined as $P =P_1 \times \ldots P_N \times \R^+$ where $P_i$ denotes the phase space of the $i$-th network node. Furthermore, let $\pi_i : P \rightarrow P_i$ be projections of (\ref{eqn:coupled_sys}), then it must hold that $\pi_i(X(t)) = x_i(t)$.
\end{Definition}

\subsection{Analysis and control}

Let: $P_1 \subseteq \R^{n_1}$, $P_2 \subseteq \R^{n_2}$, $\ldots $,  $P_N \subseteq \R^{n_N}$ be convex subsets, $P=P_1\times \ldots \times P_N$, $X=[x_1^T,\ldots , x_N^T]^T$, $x_i \in P_i$, $\phi _i :  P \times \R^+ \rightarrow P_i$ be smooth functions. In what follows, we consider systems of the form:
\begin{equation}\label{eqn:inter_sys}
\dot x_i =  \phi_i(X,t) = f_i(x_i,t) + \tilde h_i(X,t)
\end{equation}
with $i=1,\ldots , N$. Notice that (\ref{eqn:inter_sys}) represents an interconnected system (Definition \ref{defn:interconnected_system}). Specifically, in (\ref{eqn:inter_sys}) the function $f_i : P_i\times \R^+ \rightarrow P_i$ is the intrinsic dynamics of the $i$-th node, while the function $\tilde h_i : P\times \R^+ \rightarrow P_i$ describes the interaction of the $i$-th node with the other nodes composing the interconnected system.

Notice that the above formalization allows us to consider within a unique framework directed and undirected networks, self loops and multiple interactions. We will also  consider networks with (smoothly) changing topology.

The main idea of this Section for the study of the collective poly-synchronous behavior emerging in network (\ref{eqn:inter_sys}) can be stated as follows: (i) study symmetries of (\ref{eqn:inter_sys}) to determine the possible patterns of synchrony; (ii) determine among the possible patterns, the one exhibited by (\ref{eqn:inter_sys}) using contraction properties.

Consider a partition of the $N$ nodes of a network into $k$ groups, $\mathcal{G}_{1},\ldots,\mathcal{G}_{k}$, characterized by the same intrinsic dynamics. We define the following invariant subspaces associated to each group of nodes:
$$
\sM_{p,s} = \left\{x_i=x_j, \quad \forall i,j \in \mathcal{G}_s\right\}, \quad s =1,\ldots,k
$$
Notice that all the nodes of the $i$-th group are synchronous if and only if network dynamics evolve onto the associated subspace $\sM_{p,i}$. The \emph{poly-synchronous subspace}, say $\sM_p$, is then defined as the intersection of all $\sM_{p,s}$, i.e. $\sM_p = \cap_s \sM_{p,s}$, or equivalently
$$
\sM_p = \left\{x_i= x_j, \forall i,j \in \mathcal{G}_{m}, 1\le m \le k \right\}
$$
We a say that a pattern of synchrony is possible for the network of interest if its corresponding poly-synchronous subspace is flow invariant. In this view, a useful result is the following:

\begin{Theorem} \label{thm:poly_synchr}
The set $\sM_{p}$ is invariant for network (\ref{eqn:coupled_sys}) if the nodes belonging to group $\mathcal{G}_{p}$: i) have the same uncoupled dynamics; ii)  are input-equivalent.
\end{Theorem}

The proof of the above result can be found in e.g. \cite{Gol_Ste_03,Gol_Ste_06}. In terms of network synchronization, intuitively such a result implies that a specific pattern of synchrony is possible if  the aspiring synchronous nodes have synchronous input sets.

The following result is a straightforward consequence of the results of the previous Section on spatial symmetries. 

\begin{Corollary}
Assume that for network (\ref{eqn:inter_sys}) the sets  $\sM_{p,s}$ exist. Then, the synchrony pattern exhibited by the network is given by: (i)  $\sM_{p}$, if the network is contracting, or contracting towards each $\sM_{p,s}$; (ii) $\sM_{p,s}$, if the network is contracting towards $\sM_{p,s}$.
\end{Corollary}

\begin{figure}[thbp]
\begin{center}
  \includegraphics[width=5cm]{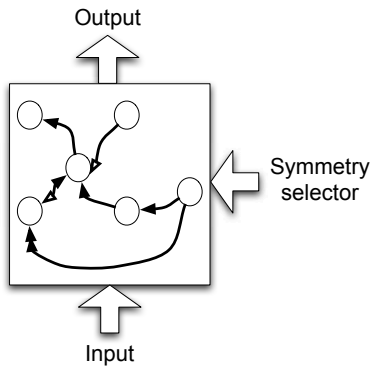}
  \caption{A schematic representation of a network used to process information. In our framework, the network is subject to two kinds of inputs: $\bar U(t)$ and $\tilde U(X)$. In particular, $\bar U(t)$ can be seen as the \emph{information} processed by the network: such input can \emph{activate} some intrinsic symmetries of the network. On the other hand, $\tilde U(X)$ is a typically sparse input that acts on the coupling functions so as to \emph{force} the activation of some desired symmetry of the interconnected system.}
  \label{fig:symmetry_selector}
  \end{center}
\end{figure}

One of the applications where the above results can be used is that of
designing networks performing specific tasks. For example, in
\cite{Pha_Slo_07} it was shown that a network with a specific symmetry
can be used to detect symmetries of e.g. images. Our results can be
used to extend this framework. Indeed, each network node (\ref{eqn:inter_sys}) may be used to process some exogenous input, $U(t) =[u_1,\ldots, u_N]$, i.e.
$$
\dot x_i =  \phi_i(X,t) = f_i(x_i,t) + \tilde h_i(X,t) + u_i(t)
$$
Now, while $u_i$ denotes the information that has to be processed by node $i$, the couplings $\tilde h_i$ may be seen as an input (typically, sparse) acting on the couplings between nodes, so as to activate a desired, arbitrary, symmetry. The
\emph{output} of the network is then some desired synchronous pattern
which arises from the intersection of the symmetries \emph{activated}
by $\bar U(t)$ and those activated by $\tilde U(X)$.  Figure
\ref{fig:symmetry_selector} schematically illustrates this principle.

For instance, assume that the system consists of a large number of
synchronized oscillators. Then we know~\cite{Wan_Slo_05} that with an
adequate choice of coupling gains:
\begin{itemize}
\item adding a single inhibitory connection between any two nodes will make the entire system contracting and therefore will stop the oscillations.
\item adding a "leader" oscillator (i.e., an oscillator with only feed-foward connections to the rest of the network, its neighbors for instance) will make the entire system get in phase with the leader.
\end{itemize}

\noindent Thus very sparse feedback inputs can completely change the symmetries of the system,
and therefore its symmetry detection specifications.

\subsubsection*{Chain topologies: revised}

Consider, again, the network topology in Figure \ref{fig:chain_topology}. Recall that in Section \ref{sec:chain_top} we proved network synchronization in two subsequent steps. Specifically, we first proved that all network trajectories are globally exponentially convergent towards the poly-synchronous subspace where $x_1=x_4$, $x_2=x_3$. We then showed that network dynamics reduced on such a subspace were globally exponentially convergent towards the synchronous subspace. 

The subspaces $\sM_2$ and $\sM_1$ towards which convergence was proved were, in turn, determined by equivariance of network dynamics with respect to some permutation action. Notice that this equivariance property is a direct consequence of the fact that node $1$ of the network is input-equivalent to node $4$ and node $2$ is input-equivalent to node $3$. Moreover, the \emph{equivalent} nodes of the $2$-nodes reduced network are also input-equivalent.

\section{Applications}\label{sec:applications}

\subsection{Synchrony patterns for distributed computing}\label{sec:hopfield}

We now turn our attention to the problem of imposing some poly-synchronous behavior for a network of interest. Specifically, we will impose different patterns of synchrony for a network composed of Hopfield models. The motivation that we have in mind here is that of \emph{multi-purpose} networks, i.e. networks that can be \emph{reused} to perform different tasks. For example, this may be the case of  sensor networks (\cite{Mar_08}, \cite{Yar_Kus_Sin_Ran_Liu_Sin_05}) where each poly-synchronous steady steady is associated to a specific set of inputs. A further notable example is the brain, where different poly-synchronous behaviors are believed to play a key role in e.g. learning processes (see e.g. \cite{Izh_06}).

We consider here a network of Hopfield models \cite{Hop_82}, \cite{Osh_Oda_07}:
\begin{equation}\label{eqn:hopfield}
\dot x_i = - x_i + \sum_{j\in N_i}a_{ij}(t)h_{ij}(x_i,x_j,t) + u_i
\end{equation}
where $a_{ij}(t)$ is the $i$-th element of the time-varying interconnection matrix $A(t)$, $h_{ij}$ represents the interconnection function  from node $j$ to node $i$ and $u_i$ in an exogenous input to the $i$-th node.

\begin{figure}[thbp]
\begin{center}
  \includegraphics[width=6cm]{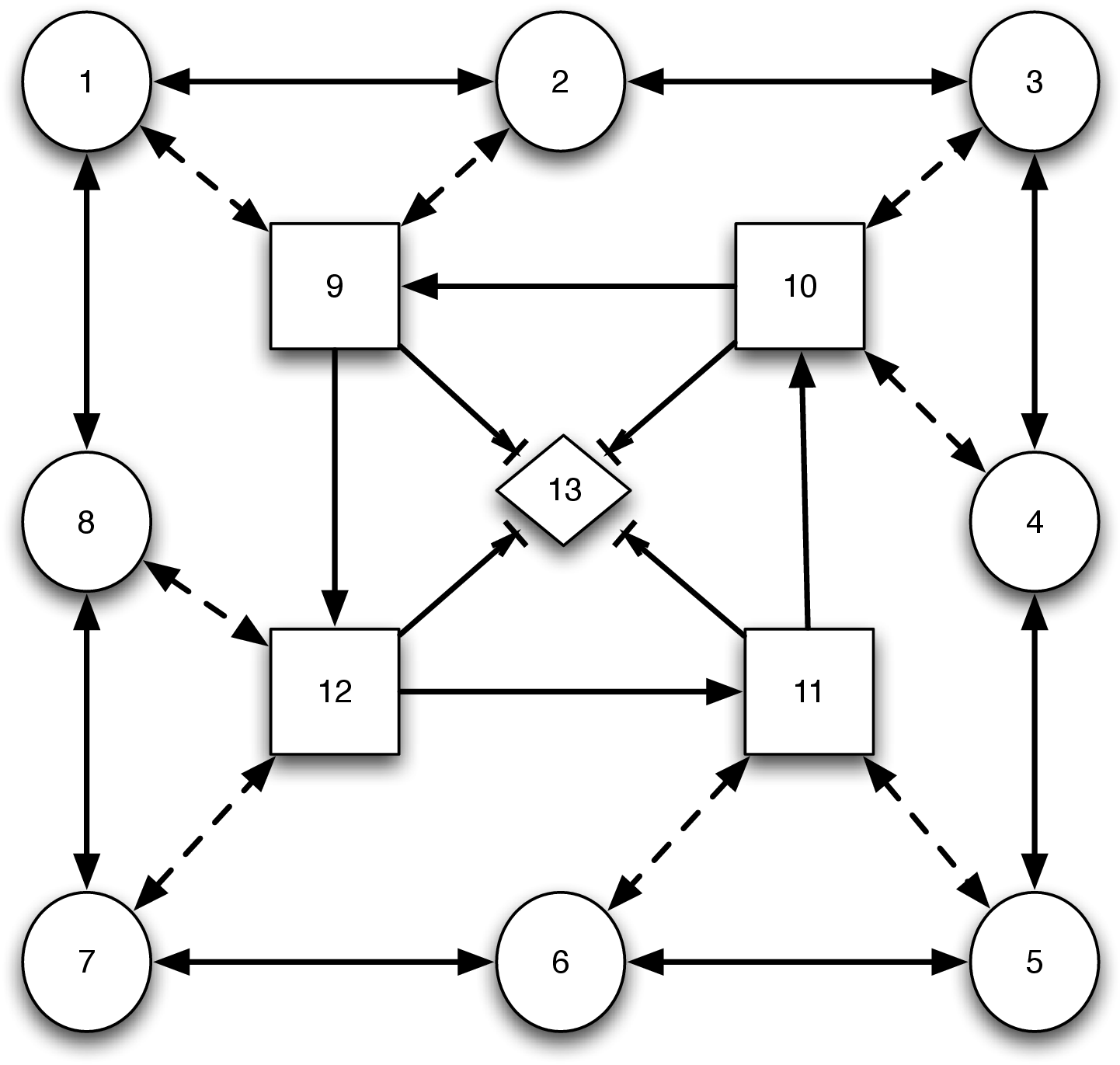}
    \includegraphics[width=6cm]{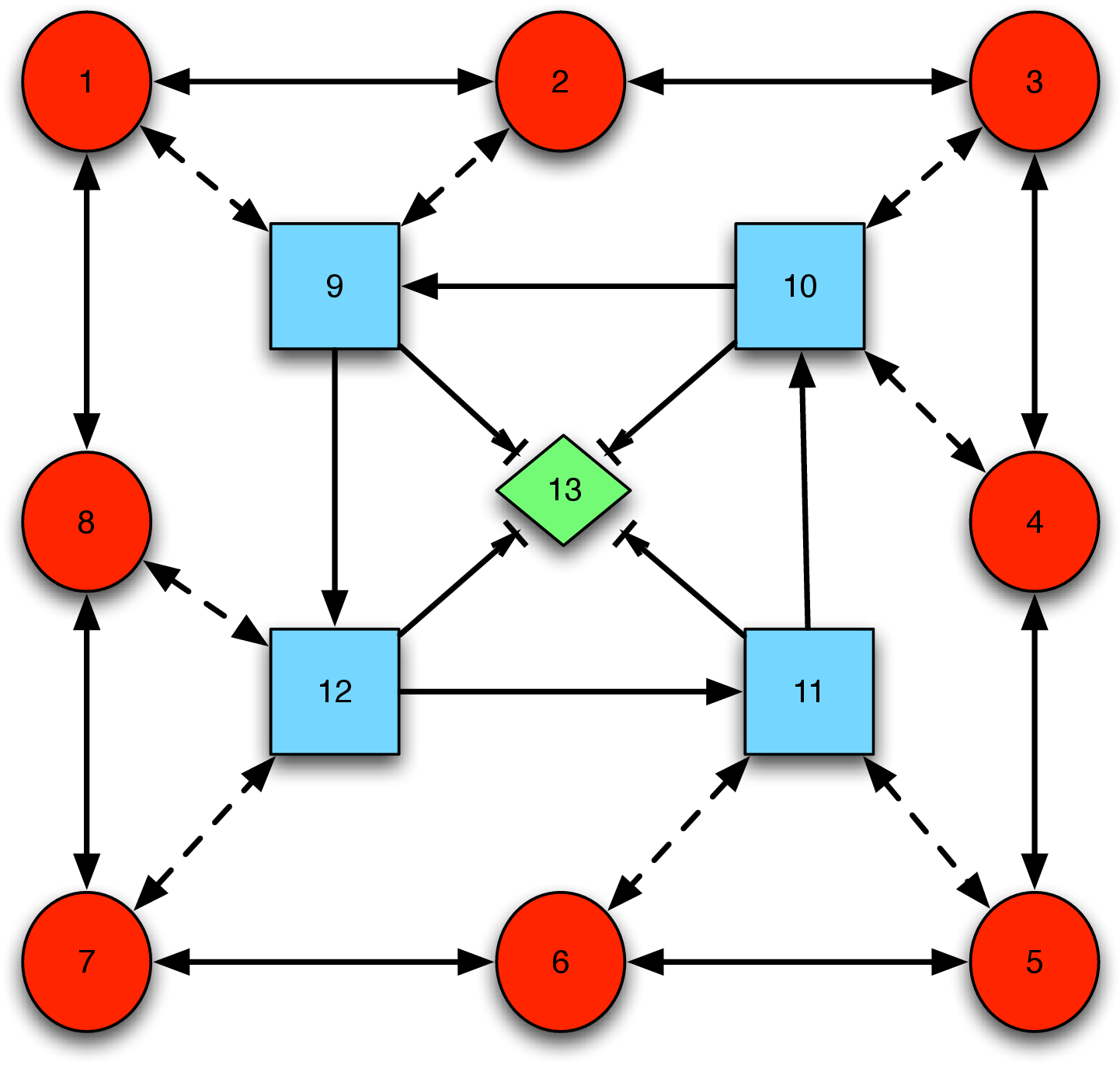}
  \caption{Network of Hopfield models used in Section \ref{sec:hopfield} (top panel). The input-equivalent nodes are pointed out in the bottom panel.}
  \label{fig:hopfield_network}
  \end{center}
\end{figure}

We start with the network in Figure \ref{fig:hopfield_network}. Nodes denoted by the same shape are forced by the same exogenous input. Specifically: (i) $u_i(t) = 1+\sin(0.7t)$ for the circle nodes; (ii) $u_i(t) = 5+3\sin(0.5t)$ for the square nodes; (iii) $u_i(t)=0$ for node $13$.

Analogously, identical arrows denote identical coupling functions:
\begin{itemize}
\item the coupling between circle nodes is diffusive, bidirectional and linear: $h_{ij}(x_i,x_j,t) = a_{ij}(t)(x_j-x_i)$
\item the coupling between square nodes is diffusive, unidirectional and linear;
\item the coupling between circle and square nodes is diffusive, bidirectional and nonlinear: $h_{ij}(x_i,x_j,t) = a_{ij}(t)\left(\arctan(x_j)-\arctan(x_j)\right)$
\item the square nodes affect the dynamics of node $13$ unidirectionally. Specifically, the dynamics of $x_{13}$ is given by:
\begin{equation}\label{eqn:coupling_central}
\dot x_{13} = -x_{13} + (1-b(t))\sum_{j=9}^{12} \frac{x_j}{1+x_j} + b(t)\sum_{j=9}^{12} \frac{1}{1+x_j}
\end{equation}
where $b(t)$ is a parameter that is smoothly increased between $0$ and $1$. Notice that $b(t)$ can be used to switch between two different coupling functions.
\end{itemize}

We remark here that the input to node $13$ is a well known coupling mechanism in the literature on neural networks, and is termed as \emph{excitatory-only} coupling, see e.g. \cite{Rub_06}.
 
It is straightforward to check that network dynamics are contracting (using e.g. the matrix measure induced by the $1$-norm).

\begin{figure}[thbp]
\centering \psfrag{x}[c]{{Time}}
\centering \psfrag{y}[c]{{$x_{i}$}}
\begin{center}
  \includegraphics[width=6cm]{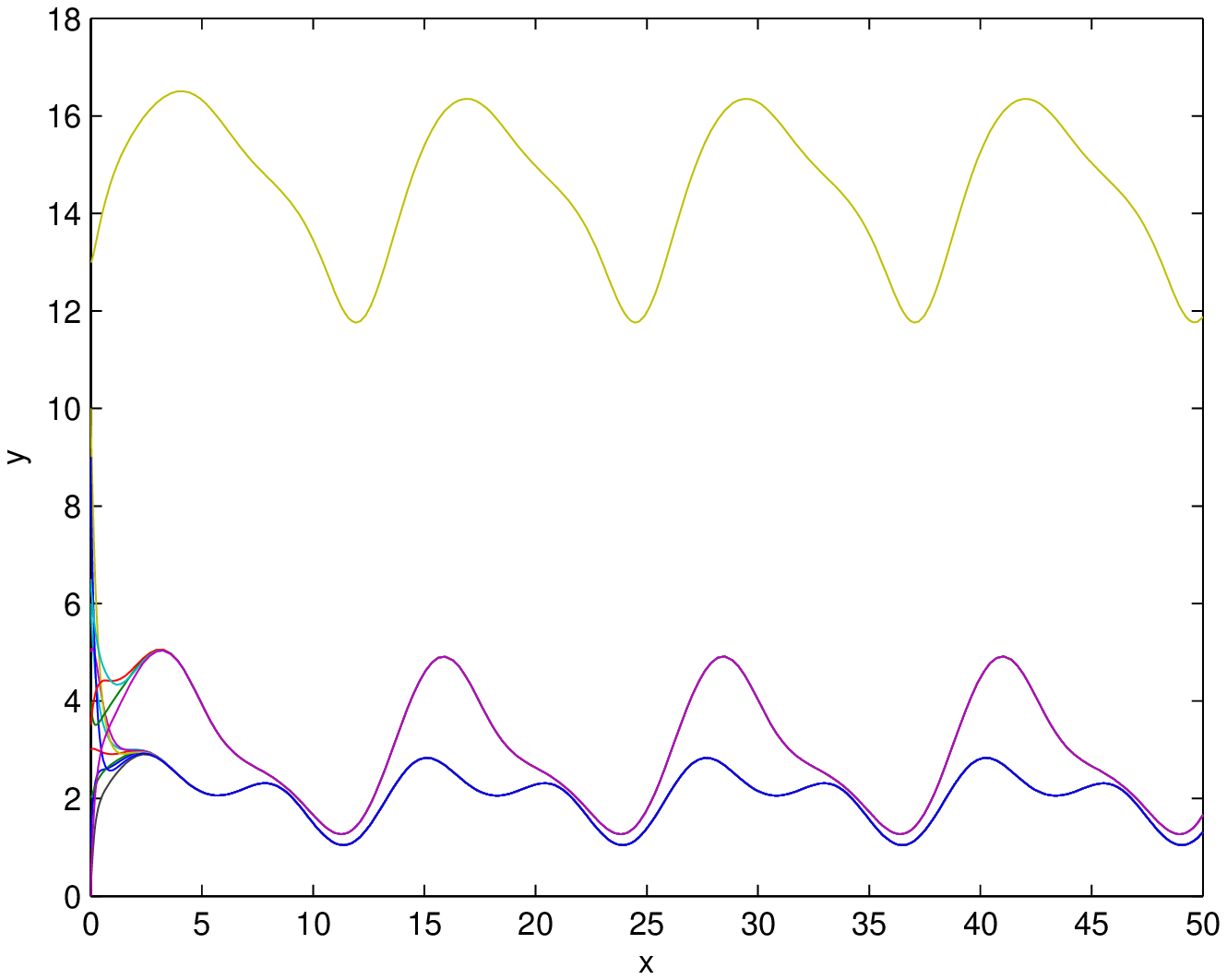}
  \centering \psfrag{x}[c]{{Time}}
\centering \psfrag{y}[c]{{$x_{i}$}}
  \includegraphics[width=6cm]{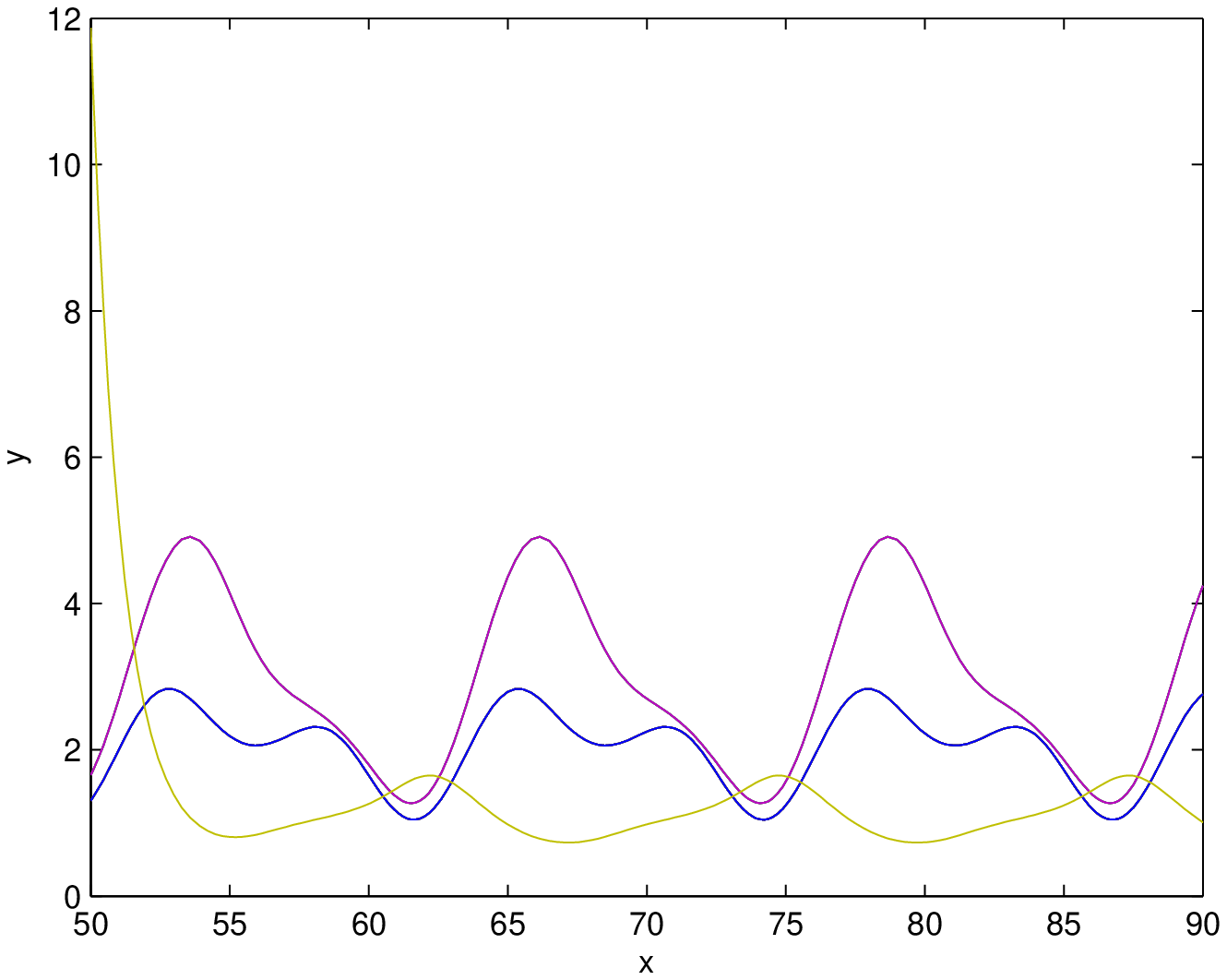}
  \caption{Network of Hopfield models (\ref{eqn:hopfield}) when (top panel) $b(t)$ in (\ref{eqn:coupling_central}) is equal to $0$. Notice the presence of three synchronized groups of nodes, corresponding the three classes of input-equivalent nodes pointed out in Figure \ref{fig:hopfield_network}. The same synchronized groups are present when $b(t)$ in is switched to $0$ (bottom panel). Notice that, in this case, the change of the coupling modifies the temporal behavior of node $13$.}
  \label{fig:simulation_hopfield_0_50}
  \end{center}
\end{figure}

In Figure \ref{fig:hopfield_network} (right panel) the input-equivalent nodes are pointed out by means of colors: the associated linear poly-synchronous subspace is
$$
\begin{array}{*{20}l}
\sM_1 = \left\{x_i=x_j, i,j=1,\ldots, 8\right\} \bigcap \\
\bigcap \left\{x_i=x_j, i,j=9,\ldots, 12\right\} 
\end{array}
$$
Furthermore, it is easy to check that $\sM_1$ is flow invariant. Now, since network dynamics are contracting, all of its trajectories converge towards a unique solution embedded into $\sM_1$. That is, at steady state all the nodes having the same \emph{color} in Figure \ref{fig:hopfield_network} are synchronized. Figure \ref{fig:simulation_hopfield_0_50} (left panel) clearly confirms the theoretical analysis, showing the presence of the three synchronized clusters, when $b(t)=0$.

The same synchronized behavior is kept even when $b(t)$ smoothly varies from $0$ to $1$. Indeed, network dynamics is still contracting and the input-equivalence propertydefining $\sM_1$ is preserved. In Figure \ref{fig:simulation_hopfield_0_50} (right panel) the behavior of the network is shown when at $t=50$, $b(t)$ is set to $1$.

Notice that the variation of $b(t)$ from $0$ to $1$ causes an \emph{inhibitory} effect of the level of $x_{13}$. This is due to the fact that, when $b(t)=0$, $x_{13}$ is forced by the sum of increasing sigmoidal functions. Vice-versa, when $b(t)=1$, $x_{13}$ is forced by the sum decreasing sigmoidal functions.

Now, assume that we need to create a new synchronized cluster consisting of e.g. nodes $2$, $4$, $6$, $8$. A way to achieve this task is that of modifying the input-equivalence property defining $\sM_1$ and to impose a new input-equivalence defining the subspace
$$
\begin{array}{*{20}l}
\sM_2 = \left\{x_i=x_j, i,j=1,3,5,7 \right\} \bigcap \\
\bigcap \left\{x_i=x_j, i,j=2,4,6,8,\right\}\bigcap\\
\bigcap \left\{x_i=x_j, i,j=9,\ldots, 12\right\}
\end{array}
$$
In turn, this can be done by smoothly varying the topology of the network, e.g. by diffusively coupling node  $13$ to the nodes $2$, $4$, $6$, $8$. The coupling function used to this aim, which preserves the contracting property, is:
$$
h_i(x_i,x_j) = h(x_j)-h(x_i),\ \ \ h(x) = \frac{1-e^{-x}}{1+e^{-x}}
$$

\begin{figure}[thbp]
\begin{center}
  \includegraphics[width=6cm]{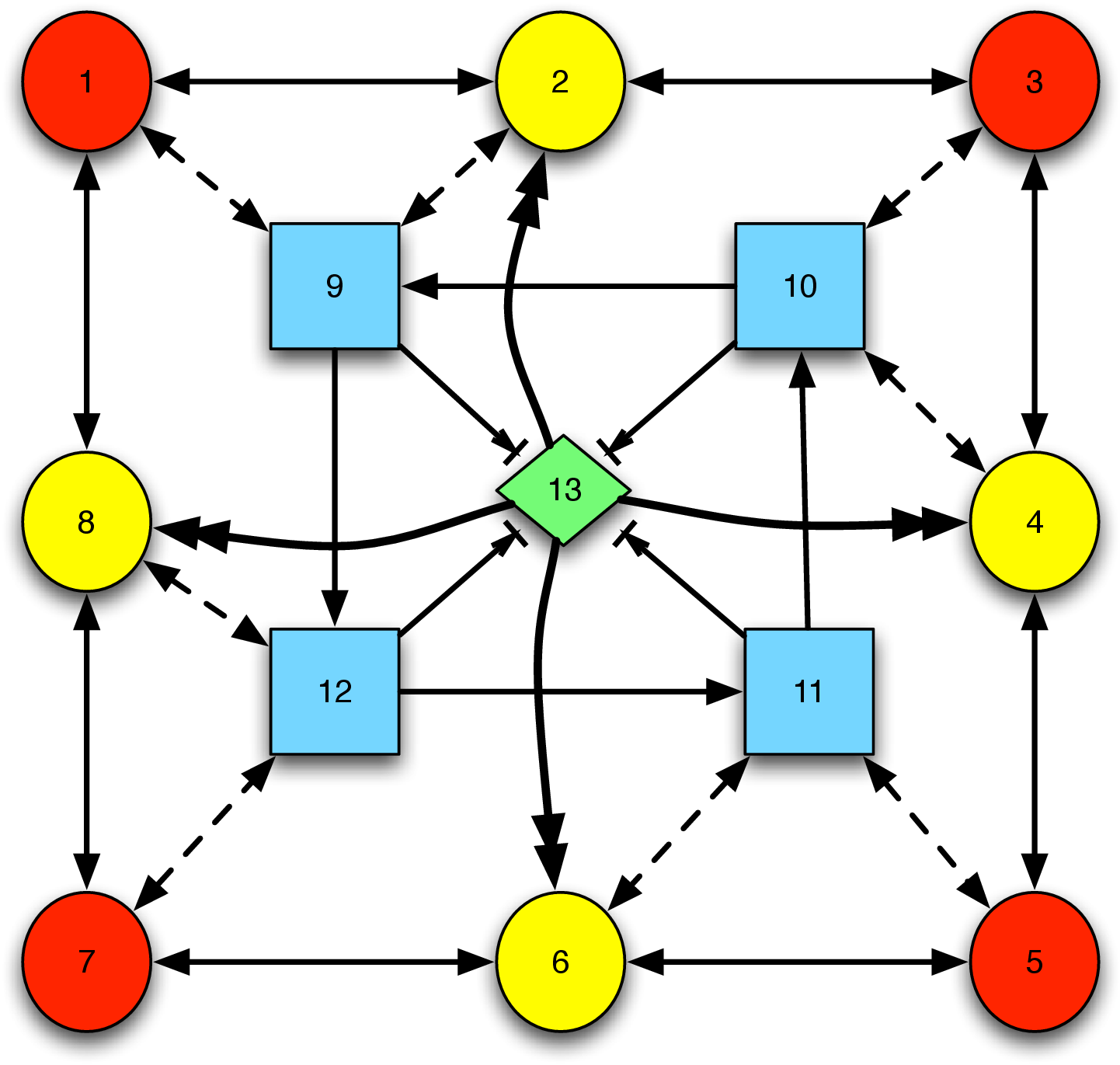}
  \centering \psfrag{x}[c]{{Time}}
\centering \psfrag{y}[c]{{$x_{i}$}}
  \includegraphics[width=6cm]{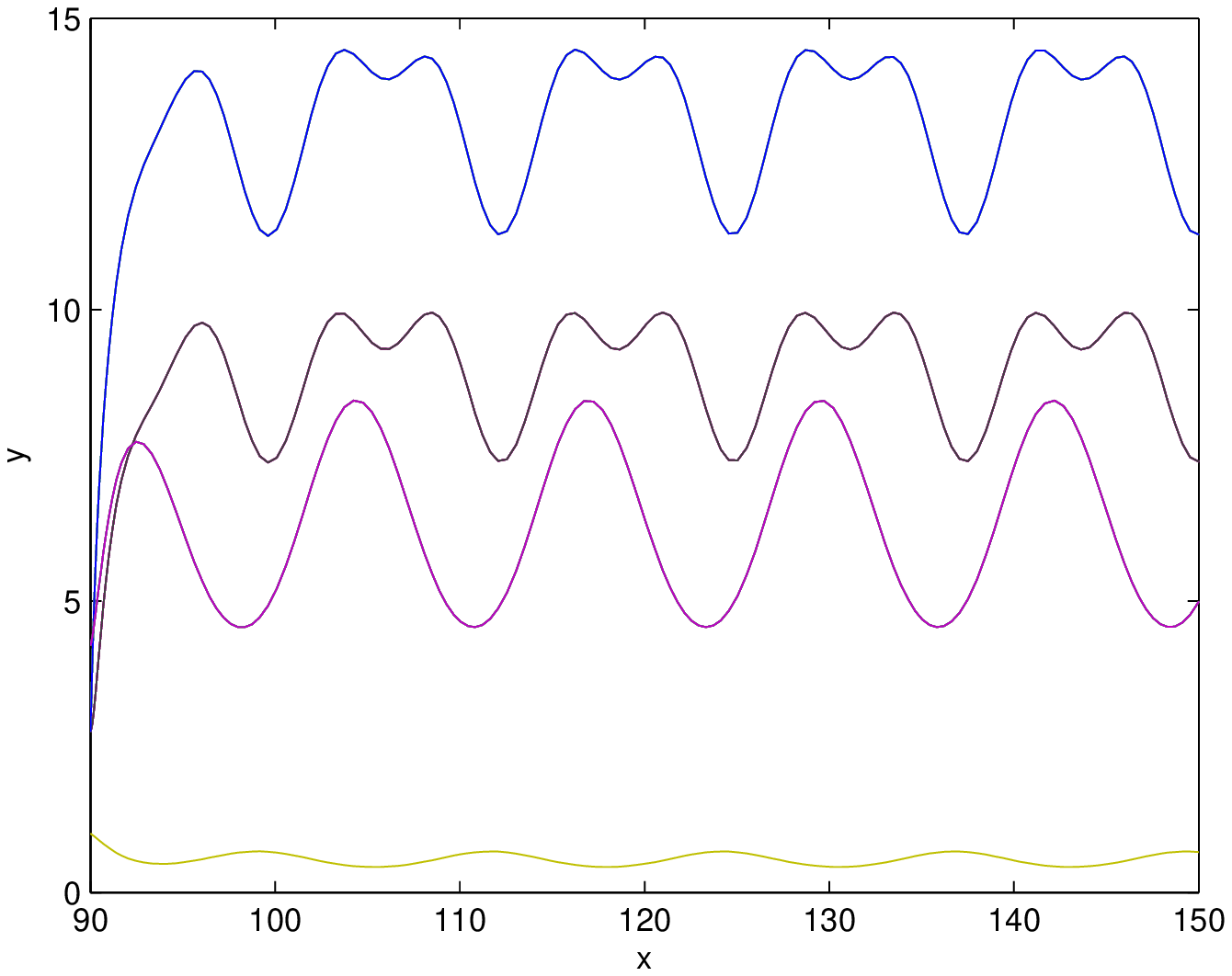}
  \caption{Network of Hopfield models used in Section \ref{sec:hopfield}. Top panel: two new links are activated by node $13$, creating a new class of input-equivalent nodes (in yellow). Bottom panel: temporal evolution of network nodes' dynamics. Notice the presence of a new group of synchronized nodes, corresponding to the new input-equivalence class.}
  \label{fig:hopfield_network_input_symmetry_2}
  \end{center}
\end{figure}

In Figure \ref{fig:hopfield_network_input_symmetry_2} (left panel) the new topology is shown, together with the class of input-equivalence. The same Figure (right panel) shows the behavior of the network, pointing out that a new cluster of synchronized nodes arises.

%\begin{figure}[thbp]
%\begin{center}
%  \includegraphics[width=8cm]{hopfield_network_input_symmetry_1}
%  \caption{Network of Hopfield models used in Section \ref{sec:hopfield}: input-equivalent nodes are pointed out by means of identical colors.}
%  \label{fig:hopfield_network_input_symmetry_1}
%  \end{center}
%\end{figure}

\subsection{Chemotaxis with quorum}\label{sec:chem_synchro}

In the above examples, we assumed that each node of the network 
communicates directly with its neighbors. This assumption on the \emph{communication mode} is often made in the literature on synchronization, see
e.g. \cite{New_2003,Boc_Lat_Mor_Cha_Hwa_06,Dan_Pal_Tsi_Has_10} and references therein.  In
many natural systems, however, network nodes do not communicate
directly, but rather by means of the 
\emph{environment}. This mechanism,
known as {\it quorum
  sensing}~\cite{Mil_Bas_01,Nar_Bas_Lev_08,Rus_Slo_10} is believed to
play a key role in bacterial infection, as well as e.g. in
bioluminescence and biofilm formation \cite{Ane_Pir_Jun_09,Nad_Xav_Lev_Fos_08}. 
Although to our knowledge this has not yet been studied experimentally, plausibly quorum sensing may also play 
a role in bacterial chemotaxis. Indeed, such a mechanism would enhance robustness of the chemotactic response~\cite{Tab_Slo_Pha_09}, with respect both to 
noise (including Brownian noise)
and to the large variations in gene expressions between individual cells~\cite{Tya_10}.

From a network dynamics viewpoint, a detailed model of such a mechanism would need to
keep track of the temporal evolution of the environmental (shared) quantity, resulting
in an additional set of ordinary differential equations \cite{Kat_08,Rus_Slo_10}:
\begin{equation}\label{eqn:general_model}
\begin{array}{*{20}l}
\dot x_ i = f\left(x_i,z\right) & i =1, \ldots, N\\
\dot z = g\left(z, \Psi\left( X\right),t\right)
\end{array}
\end{equation}
In the above equation, $N$ is the number of nodes sharing the same environment (medium). The set of state variables of the nodes is $x_i$, $X=\left[x_1^T,\ldots, x_N^T\right]^T$, while the set of the state variables of the common medium
dynamics is $z$. Notice that the medium dynamics and the medium dynamics can be of different dimensions (e.g. $x_i \in\R^n$, $z \in \R^d$).The dynamics of the nodes affect the dynamics of the common medium by means of some coupling (or input) function, $\Psi :\R^{Nn}\rightarrow \R^d$. We assume that $\partial f / \partial z$ is bounded (that is, all of its elements are bounded).

In \cite{Rus_Slo_10} it is shown that synchronization of (\ref{eqn:general_model}) is attained if the reduced order virtual system
\begin{equation}\label{eqn:virt_sys_red_ord}
\begin{array}{*{20}l}
\dot y = f\left(y,y_z\right) \\
%\dot y_z = g\left(y_z, \Psi\left( X\right),t\right)
\end{array}
\end{equation}
is contracting. Notice that the choice of such a reduced order virtual system is made possible by the fact that network equations (\ref{eqn:general_model}) are symmetric with respect to any permutation of nodes state variables, $x_i$.
  
Consider, again, the chemotaxis model (\ref{eqn:alon_second_example_1}) coupled by means of a quorum sensing mechanism, with $\phi(u/x)=u/x$:
\begin{equation}\label{eqn:network_chemo}
\begin{array}{*{20}l}
\dot x_i  = x_i (y_i-1) + h(y_i,z)\\
\epsilon \dot y_i  = \frac{u}{x_i} - y_i\\
\dot z = g(z,\Psi(Y),t)
\end{array}
\end{equation}
In the above model subscript $i$ is used to denote the state variables of the $i$-th node and $Y= [y_1^T,\ldots, y_N^T]^T$. The $i$-th node affects the dynamics of the shared variable, $z$, by means of $y_i$. Node-to-node communication is implemented by means of the input function $h$.

Notice that the presence of the coupling term and of the medium dynamics \emph{destroys} the symmetry responsible of the invariance under input scaling. We are now interested in checking under what conditions invariance under input scaling is kept for a population the chemotaxis models in (\ref{eqn:network_chemo}). We are motivated by the fact that, intuitively, bacteria go up a nutrient \emph{gradient} towards the nutrient's source, with little interest for the absolute nutrient concentration.

We model the interaction between bacteria and the environment with a \emph{dimerization} process: 
$$
h(y_i,z) = K(y_i-1)z
$$
Dimerization is a fundamental reaction in biochemical networks where two species combine to form a complex, as e.g, in the case of enzymes binding with
substrates~\cite{Sza_Ste_Per_06}. 

In what follows, we use our results to show that a possible environmental model that ensures invariance under input scaling is:
$$
\dot z = -z -1/N \sum_{i=1}^N x_iz
$$
That is, the model analyzed in the rest of this Section is:
\begin{equation}\label{eqn:sym_quorum}
\begin{array}{*{20}l}
\dot x_i = x_i(y_i-1) + K(y_i-1)z\\
\epsilon \dot y_i = \phi(u/x_i) - y_i\\
\dot z = -z -1/N \sum_{i=1}^N x_iz \\ 
\end{array}
\end{equation}
We show that such a network preserves invariance under input scaling, i.e. such a behavior is not lost when nodes are coupled through the medium dynamics, $z$.
 
Since our aim is to explain the onset of a symmetric synchronous behavior for the network, we will consider a virtual system which embeds the dynamics of the common medium:
\begin{equation}\label{eqn:virt_sym_quorum}
\begin{array}{*{20}l}
\dot y_x = y_x(y_{y}-1) + K(y_{y}-1)y_z\\
\epsilon \dot y_{y} = \phi(u/y_{x}) - y_{y}\\
\dot y_z = -y_z -1/N \sum_{i=1}^N( x_iy_z) \\
\end{array}
\end{equation}
Notice that the above system represents an hierarchy (see Section \ref{sec:math_tools}) consisting of two subsystems: the (virtual) nodes dynamics $[y_{x}(t), y_{y}(t)]$, and the (virtual) medium dynamics, $z(t)$. Therefore, (\ref{eqn:virt_sym_quorum}) is contracting if each of the subsystems is contracting (even in different metrics). 

Now, it is straightforward to check that the dynamic of $y_z$ is contracting.  Thus, to complete the analysis, we have to check if the $y$-dynamics is  contracting. Similarly to the analysis of Section \ref{sec:chemo}, we have:
\begin{equation}\label{eqn:node_dyn_rec}
\ddot y_{y}^{\ast}  + \frac{\dot y_{y}^{\ast}}{\epsilon} + \frac{1}{\epsilon}\left(\frac{u}{y_{x}} + \frac{uKy_z}{y_{x}^2}\right) y^{\ast}_{y} = 0
\end{equation}
where $y^\ast_{y} = y_{y}-1$. That is, following exactly the same steps as those used in Section \ref{sec:chemo}, we have that the above dynamics is contracting if
\begin{equation}\label{eqn:chemo_synchro_cond}
\frac{u}{y_{x}} + \frac{uKy_z}{y_{x}^2} < \frac{1}{2\epsilon}
\end{equation}
Notice that the above condition is satisfied if: i) the concentration of $y_{x}$ (and hence of $x_i$'s) is sufficiently high with respect to the concentration of $y_z$ (and hence of $z$); ii) its dynamics is sufficiently slow ($\epsilon$ small); iii) $K$ is properly tuned. Recall from Section \ref{sec:chemo}  that both i) and ii) are true by hypotheses.

Thus, under the above conditions we have that the virtual system is contracting, and hence synchronization of network nodes is attained. Moreover, it is straightforward to check that the particular choice of the coupling function and of the medium dynamics ensures that all the hypotheses of Theorem \ref{thm:virtual_system_diff_symmetries} are satisfied with the actions
\begin{equation}\label{eqn:symmetries_quorum}
\begin{array}{*{20}l}
\gamma_i (y_{x,i}, y_{y,i}, y_z) \rightarrow (y_{x,i}/u_{\min,i}, y_{y,i}, y_z/u_{\min,i})\\
\rho_i u_i \rightarrow u_i/u_{\min,i}\\
\end{array}
\end{equation}
where the subscripts $i$ denote the system state variables and the actions associated to the input $u_i(t)$ (see Section \ref{sec:chemo} for further details). This, in turn implies that the virtual system exhibits the invariance under input scaling. Now, recall that network nodes are particular solutions of the virtual system. Thus, all network nodes globally exponentially converge towards each other (since the virtual system is contracting), while exhibiting the invariance under input scaling (by virtue of Theorem \ref{thm:virtual_system_diff_symmetries}).

\subsection{Quorum with periodic inputs or communication delays}

Quorum sensing mechanisms exploit the symmetry of a dynamic system
under permutation of individual elements. From a contraction point of
view, this allows one to use a virtual system of the same dimension as
individual elements, and such that each individual trajectory
represents a particular solution of the virtual system.

This principle extends straightforwardly to two cases of practical
interest. The first is the case when the environment or the entire
system is also subject to an external periodic input, thus yielding a
spatio-temporal symmetry.  The second case is when communications
between nodes and the environment exhibit significant delays. This
case may model actual delays information transmission (e.g., in a
natural or robotic swarm application) and signal processing, or for
instance the effect of diffusion or of non-homogeneous concentrations
in a biochemical context.

We now show that the combined use of symmetries and contraction analysis can be used to provide a sufficient condition to control the periodicity of the synchronous final behavior of a quorum sensing of interest. The idea is to force a network of interest by means of a periodic input, $r(t)$, and then to provide conditions ensuring that the network becomes synchronized onto a periodic orbit having the same period as $r(t)$. See the companion paper \cite{Rus_Slo_10} for further details. Our main result is as follows.
\begin{Theorem}\label{thm:control_basic}
Consider the following network
\begin{equation}\label{eqn:general_mode_controll}
\begin{array}{*{20}l}
{\dot x}_i = \tilde f\left(x_i,z\right) & \ \ \ \ \ i =1, \ldots, N\\
\dot z = \tilde g\left(z,\Psi\left(x_1,\ldots,x_N\right)\right) + r\left(t\right)
\end{array}
\end{equation}
where $r\left(t\right)$ is a $T$-periodic signal.  All the nodes of the network synchronize onto a periodic orbit of period  $T$, say $x_T(t)$, if: (i) $f\left(x_i,v(t)\right)$ is a contracting function; (ii) the reduced order system ($x_c(t) \in \R^n$)
$$
\begin{array}{*{20}l}
{\dot x_c} = f\left(x_c,z\right) \\
\dot z = g\left(z,\Psi\left(x_c,\ldots,x_c\right)\right) + r\left(t\right)
\end{array}
$$
is contracting.
\end{Theorem}
\noindent 
Note that the dynamics $\tilde f$ and $\tilde g$ include the coupling terms between nodes and environment.
\proof
Consider the virtual system
\begin{equation}\label{eqn:virtual_quorum_controlled}
\begin{array}{*{20}l}
\dot y_1 = f\left(y_1,z\right)\\
\end{array}
\end{equation}
By hypotheses (\ref{eqn:virtual_quorum_controlled}) is contracting and hence the nodes state variables will converge towards each other, i.e. $\abs{x_i -x_j}\rightarrow 0$ as $t\rightarrow + \infty$. That is, all the network trajectories converge towards a unique common solution, say $x_c(t)$. This in turn implies that, after transient, network dynamics are described by the reduced order system
$$
\begin{array}{*{20}l}
{\dot x}_c = f\left(x_c,z\right) \\
\dot z = g\left(z,\Psi\left(x_c,\ldots,x_c\right)\right) + r\left(t\right)
\end{array}
$$
Now, the above system is contracting by hypotheses and $r(t)$ is a $T$-periodic signal. In turn, this implies that all of its solutions will converge towards a unique $T$-periodic solution, i.e.
$$
\abs{x_c(t) - x_T(t)} \rightarrow 0, \ \ \ \ \ t \rightarrow + \infty
$$
This proves the result.
\endproof
Notice that the proof of the above result is based on the combined use of symmetries and contraction. Specifically, the use of a reduced order virtual system is made possible by the fact that the network is symmetric with respect to any permutation of the nodes state variables. Moreover, contraction analysis provides sufficient conditions guaranteeing that all the solutions of the virtual system converge to a periodic trajectory having the same period as the input, $r(t)$. Since the nodes' state variables are particular solutions of the virtual system, this implies that network nodes are synchronized onto a periodic orbit having the same period of $r(t)$.

The case when some form of delay occurs in the communication can be
treated similarly, by using results on the effect of delays in
contracting systems~\cite{Wan_Slo_06}.
Consider for instance a network of linearly diffusively coupled nodes $x_i$ coupled by means of a quorum sensing mechanism,
\begin{equation}\label{eqn:quorum_delay}                                                                      
\begin{array}{*{20}l}                                                                                         
\dot x_i = f(x_i) +  K_{iz} (z(t-T_{zx}) - x_i(t) ) \\                                                                  
\dot z = g(z)+ \frac{1}{N} \sum K_{zi}(x_i(t-T_{xz}) - z(t) ) \\                                                    
\end{array}                                                                                                   
\end{equation}
where: $f(\cdot)$ and $g(\cdot)$ (denoting the intrinsic dynamics of the network nodes and of the common medium) are contracting within the same metric, the constant $T_{zx} \ge 0 $ represents a communication or computation delay from the medium to the nodes, and similarly the constant $T_{xz} \ge 0$ represents a delay from the nodes to the medium.

Notice that network (\ref{eqn:quorum_delay}) is symmetric with respect to any permutation of the nodes' state variables. As discussed above, this symmetry implies that the network can be analyzed by using the reduced order virtual system
\begin{equation}\label{eqn:quorum_delay_virtual}                                                                      
\begin{array}{*{20}l}                                                                                         
\dot y_x = f(y_x) +  K_{iz} (y_z(t-T_{zx}) - y_x(t) ) \\                                                                  
\dot y_z = g(y_z)+ \frac{1}{N} \sum K_{zi}(y_x(t-T_{xz}) - y_z(t) ) \\                                                    
\end{array}                                                                                                   
\end{equation}
As proven in \cite{Wan_Slo_06}, all the trajectories of the above system converge towards each other if $f(\cdot)$ and $g(\cdot)$ are contracting within the same metric. Since the nodes' state variables are particular solutions of the virtual system (\ref{eqn:quorum_delay_virtual}), it follows that all the solutions of the network converge towards a fixed point in the network phase space. Now, as shown in \cite{Wan_Slo_06}, if $K_{iz} = \frac{1}{N}K_{zi}$, then $z(t)$ tends to the fixed value $\bar z$ and all nodes tend to the common equilibrium value $ x_i = x_j = \bar x \ ,\ \ \forall i,j$. Moreover, $\bar x$ and $\bar y$ are such that:
$$
\begin{array}{*{20}l}
f(\bar x) + k_{iz} \left(\bar z - \bar x\right) = 0\\
g(\bar z) + k_{iz} \left(\bar x - \bar z\right) = 0\\
\end{array}
$$ 

\section{Concluding remarks}\label{sec:conclusions}

We presented a framework for analyzing/controlling stability of
dynamical systems by using a combination of structural properties of
the system's vector field (symmetries) and convergence properties
(contraction). We first showed that the symmetries of a contracting
vector field are transferred to its solutions, and then generalized
this result by using virtual systems. In turn, those results were used
to describe invariance under input scaling in transcriptional systems,
a property believed to play a key role in many sensory systems. The
case of quorum-sensing network with delays was also considered.
Finally, we showed how our results could suggest a mechanism for
quorum sensing in bacterial chemotaxis, thus combining symmetries in
cell interactions with invariance to input scaling.

%\bibliography{russobib}

%merlin.mbs 2010-03-15 4.21a (PWD, AO, DPC)

 \end{document}